\documentclass{elsarticle}
\usepackage{graphics,xcolor}
\usepackage[normalem]{ulem}
\usepackage{amssymb}
\usepackage{amsmath}
\usepackage{amsfonts}
\usepackage{epsfig,color}
\usepackage{tikz}
\usepackage[mathscr]{euscript}
\usepackage{epstopdf}
\usepackage{ntheorem}
\usepackage{caption}
\usepackage{subcaption}
\usepackage{booktabs}
\usepackage{wrapfig}
\usepackage{hyperref}
\usepackage{array,multirow}

\font\escript=eusm10

\newenvironment{proof}{\paragraph{Proof}}{\hfill$\square$}

\newtheorem{theorem}{Theorem}[section]
\newtheorem{lemma}[theorem]{Lemma}

\newtheorem{remark}[theorem]{Remark}

\newtheorem{example}[theorem]{Example}

\catcode`@=11 \@addtoreset{equation}{section}
\renewcommand\theequation{\thesection.\@arabic\c@equation}
\catcode`@=12

\makeatother

\journal{Mathematics and Computers in Simulation}

\begin{document}
	
\begin{frontmatter}
\title{Anisotropic Error Analysis of Weak Galerkin finite element method for Singularly Perturbed Biharmonic Problems}

\author[address1]{Aayushman Raina}
\ead{aayushman.raina@iitg.ac.in}
\author[address1]{Srinivasan Natesan \corref{correspondingauthor}}
\ead{natesan@iitg.ac.in}
\cortext[correspondingauthor]{Corresponding author}	
\address[address1]{Department of Mathematics, Indian Institute of Technology Guwahati, Guwahati, 781 039 - {\sc India}}
\author[address2]{{\c{S}}uayip Toprakseven}
\ead{topraksp@artvin.edu.tr}
\address[address2]{Department of Mathematics, Faculty of Art and Science, Yozgat Bozok University, 66100 Yozgat,  - {\sc Turkey}}

\begin{abstract}
We consider the Weak Galerkin finite element approximation of the Singularly Perturbed Biharmonic elliptic problem on a unit square domain with clamped boundary conditions. Shishkin mesh is used for domain discretization as the solution exhibits boundary layers near the domain boundary. Error estimates in the equivalent $H^{2}-$ norm have been established and the uniform convergence of the proposed method has been proved. Numerical examples are presented corroborating our theoretical findings.
\end{abstract}

\end{frontmatter}

\goodbreak\noindent
\section{Introduction}\label{ch2_secint}

In this paper, we are concerned to solve the following biharmonic singularly perturbed problem
\begin{equation}\label{ch1_md_problem}
\left\{\begin{array}{ll}
\varepsilon^{2}\Delta^{2} u - \Delta u + au = \textsl{g}, \,\, \mbox{in} \,\, \Omega,\\ [6pt]
u = 0, \dfrac{\partial u}{\partial n} = 0 \,\, \mbox{on} \,\, \partial\Omega,
\end{array}\right.
\end{equation}
where {$\Omega\subset \mathbb{R}^2$}  is convex, bounded and polygonal domain. Here ${0<\varepsilon\ll 1}$ is a singular perturbation parameter, and {the function} $a$  satisfies {$0\leq a_1<a(x,y)<a_2$} for some  constants $a_{1}$ and $a_{2}$. {We also} assume that the functions $\textsl{g}(x,y)$ {and $a(x,y)$ are sufficiently smooth on $\Omega$}. The problem \eqref{ch1_md_problem} models different types of physical models, {for example}, in the elasticity theory \eqref{ch1_md_problem} models a plate bending problem (clamped plate), where the parameter $\varepsilon$ measures the ratio of bending rigidity to tensile stiffness in the plate. In fluid dynamics it can be obtained from stream-vorticity formulation of Oseen equations, where $\varepsilon$ is the reciprocal of the Reynolds number. In 3D, \eqref{ch1_md_problem} can be seen as a simplified form of the stationary Cahn-Hiliard equations, with $\varepsilon$ related to the length of the transition region of phase separation.

Owing to the significance of such fourth order problems, a lot of methods have been introduced to solve such problems. Some conforming finite element methods were studied in \cite{clough1965finite,dupont1979family} which requires $\mbox{\escript{C}}^{1}$ function space of piecewise polynomials. To obtain $\mbox{\escript{C}}^{1}$ elements, high-degree polynomials are employed, which leads to an increase in computational cost and time. To overcome these challenges by avoiding $\mbox{\escript{C}}^{1}$ elements, non-conforming and mixed finite element methods were discussed in \cite{monk1987mixed,lascaux1975some}. Discontinuous Galerkin finite element methods have also been employed for solving the biharmonic problems, e.g., in \cite{mozolevski2007hp} an $hp$- version of the interior penalty discontinuous Galerkin was used and in \cite{engel2002continuous, franz2014c0} a continuous interior penalty method ($\mbox{\escript{C}}^{0}$-IPM) was developed for the biharmonic problems.

Wang and Ye first proposed a new class of finite element methods, called the weak Galerkin finite element method (WG-FEM) in \cite{wang2013weak} for a second-order elliptic problem. By introducing a concept of \textit{weak functions and weak gradients}, the WG-FEM provides a flexibility in constructing higher order elements, as one can use discontinuous basis functions on a general partition of the domain. Due to its flexibility on choosing mesh division and approximation spaces, the Weak Galerkin method has been studied and implemented on different kind of problems arising in physics and engineering. Some work on solving biharmonic problem using WG-FEM has been done in \cite{mu2014weak,ye2020new}.

Recently, many researchers have been working on implementing WG-FEM to solve singular perturbation problems (SPPs). Lin et al. in \cite{lin2018weak} considered a singularly perturbed second-order convection-diffusion-reaction problem in the convection dominated regime and solved it using WG-FEM. Zhu and Xie in \cite{zhu2020uniformly} implemented a WG-FEM for one-dimensional convection-diffusion problems on a Shishkin mesh. A $P0-P0$ WG-FEM for a diffusion-reaction SPP was discussed in \cite{al2020p0}. Toprakseven in \cite{toprakseven2022optimal} proposed a WG-FEM for solving a nonlinear singular perturbation of reaction-diffusion type on a Bakhvalov-type mesh. Optimal-order uniform convergence of the method in the energy norm and a stronger balanced norm was established.

For solving singularly perturbed biharmonic problems using WG-FEM, not much work has been done. In \cite{cui2020uniform} they have considered \eqref{ch1_md_problem} without the reaction term and have constructed one Weak Galerkin finite element which converges uniformly with respect to the perturbation parameter $\varepsilon \ll 1$. Under certain regularity assumptions , they showed optimal order of convergence in the discrete $H^{2}-$norm and in the $L^{2}-$norm, and their error estimates were based on uniform meshes. But the uniform convergence of the WG-FEM on layer-adapted mesh has not been discussed so far for fourth order problems. Our goal is to investigate the uniform-convergence of the WG-FEM on a Shishkin mesh for problem of type \eqref{ch1_md_problem}.

As SPPs show boundary layer character in its solution, standard numerical methods fail to catch the solution profile correctly, are inefficient, and produce non-physical oscillations on coarse mesh. That is why analysis and implementation of such numerical methods on layer-adapted meshes are important. In this work, for the first time, a WG-FEM has been implemented over a layer-adapted mesh for a fourth-order SPP. Anisotropic error estimates have also been derived in the equivalent norm of $H^2$ and the uniform convergence of the proposed method has also been proved.

The remainder of this paper is organized as follows:  In Section \ref{model_prob}, we present our model fourth-order singularly perturbed problem, its solution decomposition and discusses its WG-FEM formulation. Further, the piecewise uniform Shishkin mesh is constructed in Section \ref{model_prob}. In Section \ref{WG_Discre} computation of local discretizations has been done. We have obtained the corresponding error equation and the anisotropic error estimates in the equivalent norm $H^{2}$ using an interpolation operator and the projections $L^{2}$ in Section \ref{Error_Analysis}. We have presented some numerical examples to validate our theoretical findings in Section \ref{Numerical_Examples}. {In Section \ref{Conclusions}, we close this paper with conclusions. }

\section{The model problem }\label{model_prob}
\subsection{Solution decomposition}

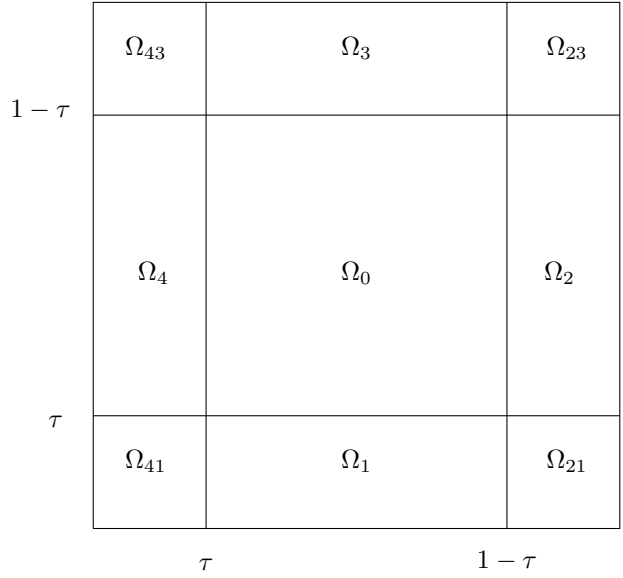
\begin{figure}
\centering
\begin{tikzpicture}
\draw (0,0) -- (7,0) -- (7,7) -- (0,7) -- (0,0);
\draw (1.5,0) -- (1.5,7);
\draw (5.5,0) -- (5.5,7);
\draw (0,1.5) -- (7,1.5);
\draw (0,5.5) -- (7,5.5);
\node[label=$\Omega_{41}$]  at (0.7,0.5) {};
\node[label=$\Omega_{21}$]  at (6.3,0.5) {};
\node[label=$\Omega_{23}$]  at (6.3,6) {};
\node[label=$\Omega_{43}$]  at (0.7,6) {};
\node[label=$\Omega_{0}$]  at (3.5,3) {};
\node[label=$\Omega_{1}$]  at (3.5,0.5) {};
\node[label=$\Omega_{2}$]  at (6.2,3) {};
\node[label=$\Omega_{3}$]  at (3.5,6) {};
\node[label=$\Omega_{4}$]  at (0.8,3) {};
\node[label=\textbf{$\tau$}]  at (1.5,-0.8) {};
\node[label=\textbf{$1-\tau$}]  at (5.5,-0.8) {};
\node[label=\textbf{$\tau$}]  at (-0.5,1.1) {};
\node[label=\textbf{$1-\tau$}]  at (-0.7,5.2) {};
\end{tikzpicture}
\caption{Discretization of unit square $\Bar{\Omega}$} \label{ch1_mesh}
\end{figure}

Consider the model problem \eqref{ch1_md_problem}.
In \cite{franz2014c0}, authors made reasonable assumptions on the decomposition of solution and the presence of boundary layers over a unit square {\em i.e.}, $\Omega = (0,1)\times (0,1)$. The solution decomposition result has been stated in the following Lemma:

\begin{lemma}[\cite{franz2014c0}] \label{solutiondecom}
The solution $u$ of the problem \eqref{ch1_md_problem} in a unit square has the following decomposition:
\[
u = \mathscr{S} + \sum_{i \in K}^{}\mathscr{E}_{i}
\]
where $K := \{1,2,3,4,21,41,23,43\}$, $\mathscr{S}$ is the smooth term and $\mathscr{E}_{i}$ are the layer terms with $\mathscr{S}, \mathscr{E}_{i} \in H^{s}(\Omega), i \in K$ for any natural number $s$ that satisfies $3 \leq s \leq k+1$. One can refer Figure \ref{ch1_mesh} for better understanding. Furthermore, for each $\tau \in (0,1/4)$, we have to follow the bounds:
\begin{eqnarray}\label{ch1_Dsol1}
|\mathscr{S}|\strut_{H^{s}(G)} \leq \mbox{meas}\,(G)^{1/2}, G \subseteq \Omega,
\end{eqnarray}
\begin{eqnarray}\label{ch1_Dsol2}
\|\mathscr{E}_{i}\|\strut_{\infty, \Omega \setminus \tilde{\Omega}_{i}} \leq C\,\varepsilon \,\exp\left(-\dfrac{\tau}{\varepsilon}\right),
\end{eqnarray}
\begin{eqnarray}\label{ch1_Dsol3}
|\mathscr{E}_{i}|\strut_{H^{s}(\tilde{\Omega}_{i})} \leq C \, \varepsilon^{\frac{3}{2}-s},
\end{eqnarray}
\begin{eqnarray}\label{ch1_Dsol4}
|\mathscr{E}_{i}|\strut_{H^{l}(\Omega \setminus \tilde{\Omega}_{i})} \leq C \, \varepsilon^{\frac{3}{2}-l}\, \exp\left(-\dfrac{\tau}{\varepsilon}\right), \quad  l \in \{1,2,3\},
\end{eqnarray}
where $\widetilde{\Omega}_{i} := \Omega_{i}, i \in \{41,21,23,43\}$.
\end{lemma}

 %%%%%%%%%%%%%%%%%%%%%%%%%%%%%%%%%%%%%%%%%%%%%%%%%%%%%%%%%%%%%%%%%%%%%%%%%%%%

\subsection{The Shishkin {triangular}  mesh}

Our estimates \eqref{ch1_Dsol1}, \eqref{ch1_Dsol2}, \eqref{ch1_Dsol3} and \eqref{ch1_Dsol4} allow us to take the same mesh transition parameter $\tau:= \min\{ 1/4, \varepsilon \lambda \ln N\}$ for both $x-$ axis and $y-$ axis, where $\lambda$ will be determined later. We consider our mesh to be the Shishkin mesh, which can be thought of as a tensor product of a one-dimensional Shishkin mesh with itself. Let $\tau$ and $1-\tau$ be our mesh transition parameters, and
\begin{eqnarray} \label{ch1_h1_h2}
h_{1}:= 4\tau/N \,\,\mbox{and}\,\, h_{2}:= 2(1-2\tau)/N,
\end{eqnarray}
where $N$ denotes the number of elements along the $x$ and $y$ axes and is assumed to be a multiple of 4.
Then we can define our mesh points as
\begin{equation}\label{ch1_xmesh}
x_{i} =
\left\{\begin{array}{ll}
ih_{1}, \quad &i = 0,\ldots,\frac{N}{4},\\ [6pt]
\tau + (i-\frac{N}{4})h_{2}, \quad &i = \frac{N}{4}+1,\ldots,\frac{3N}{4}, \\ [6pt]
1-\tau + (i-\frac{3N}{4})h_{1}, \quad &i = \frac{3N}{4} + 1,\ldots,N.
\end{array}\right.
\end{equation}
for $x$-coordinate axis. In the same fashion, on $y$-coordinate axis, we have
\begin{equation}\label{ch1_ymesh}
y_{j} =
\left\{\begin{array}{ll}
jh_{1}, \quad &j = 0,\ldots,\frac{N}{4},\\ [6pt]
\tau + (j-\frac{N}{4})h_{2}, \quad &j = \frac{N}{4}+1,\ldots,\frac{3N}{4}, \\ [6pt]
1-\tau + (j-\frac{3N}{4})h_{1}, \quad& j = \frac{3N}{4} + 1,\ldots,N.
\end{array}\right.
\end{equation}
Therefore, using these mesh points we can define our rectangular Shishkin mesh $T_{N}$ as a tensor product of 1D Shishkin mesh in each direction. That is,
\begin{eqnarray}\label{ch1_tensor_mesh}
{T}_N := \{T_{i,j} := [x_{i-1},x_{i}] \times [y_{j-1},y_{j}] \,\,\mbox{for}\,\, i,j = 1,\ldots, N\}.
\end{eqnarray}
Given a mesh rectangle, drawing the off-diagonal in each sub-rectangle, we get a piecewise uniform triangulation $\mathcal{T}_{N}$ which is our desired Shishkin triangular mesh. The mesh sizes $h_{x, i}=x_{i+1}-x_i$ and $h_{y, j}=y_{j+1}-y_j$ have the following order:
$$
\begin{aligned}
&  h_{x, i} =\mathcal{O} (N^{-1})=h_{y, j}, \quad i, j \in\{N / 4, N / 4+1 \ldots, 3N/4-1\}, \\
&  h_{x, i}=\mathcal{O}( \varepsilon N^{-1} \ln N)=h_{y, j}, \quad i, j \in\{0,1, \ldots, N / 4-1\}\cup \{3N/4, \ldots, N -1\}.
\end{aligned}
$$
Let $h_T$ be the mesh size of a triangle $T$ in our Shishkin triangulation $\mathcal{T}_N$.
\\
\begin{figure}[htb]
\centering
\includegraphics[height=9cm]{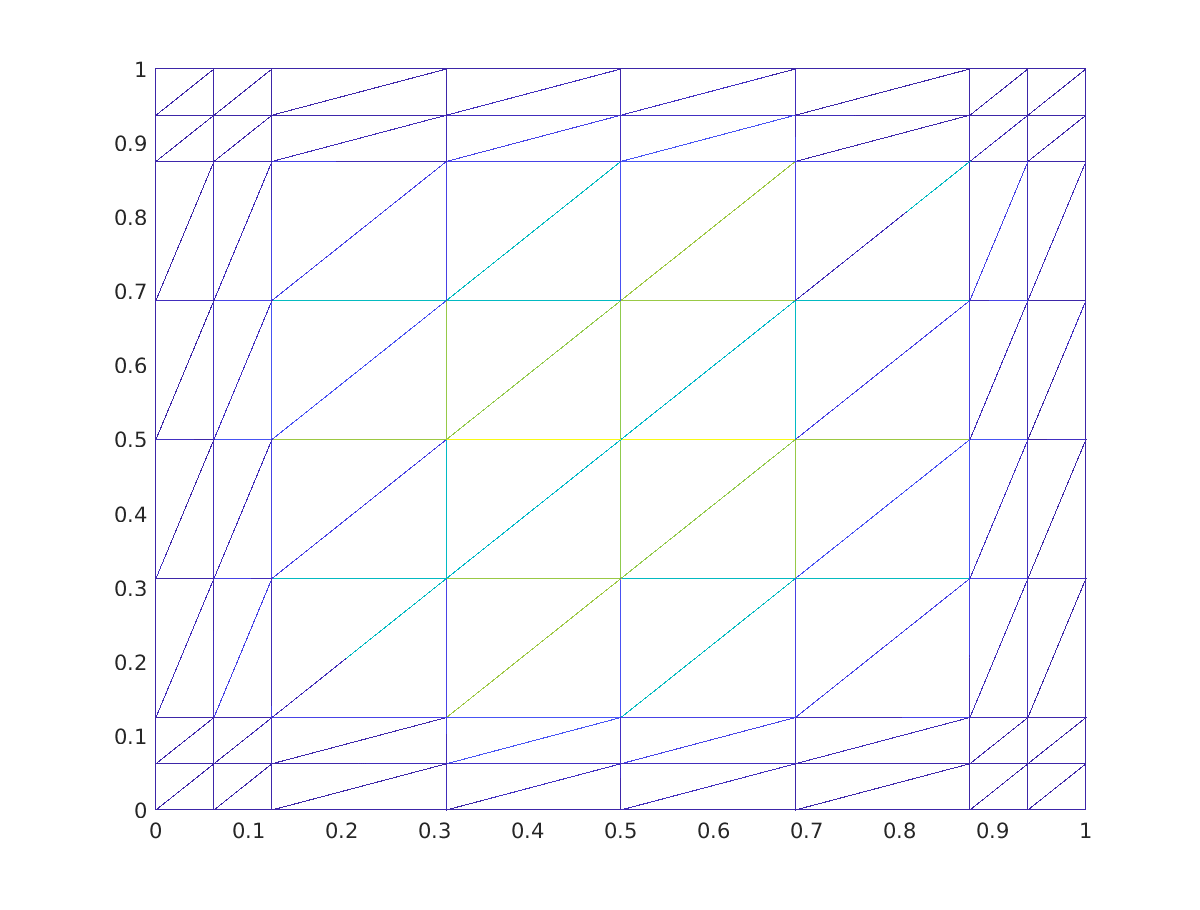}
\caption{A Shishkin triangular mesh with $N=8$}
\label{example2}
\end{figure}

\subsection{Weak formulation and well-posedness}
Consider our model problem as follows
\begin{equation*}
\left\{\begin{array}{ll}
\varepsilon^{2}\Delta^{2} u - \Delta u +au = \textsl{g}, \,\, \mbox{in} \,\, \Omega,\\[6pt]
u = 0, \quad \dfrac{\partial u}{\partial n} = 0, \,\, \mbox{on} \,\, \partial\Omega,
\end{array}\right.
\end{equation*}
where $\Omega=(0,1)\times(0,1)$, and $u \in H^{2}(\Omega) \,\,\mbox{and}\,\, \textsl{g} \in L^{2}(\Omega)$.

The WG-FEM provides high flexibility in element construction, by using discontinuous piecewise polynomials on a general discretization of the domain.
By choosing an appropriate local basis functions, the desired order of accuracy can be obtained. For fourth-order problems, WG-FEM requires the construction of the weak Laplacian operator denoted by $\Delta_{w}$ and if the convection term is also present then weak gradient operator $\nabla_{w}$ is also required.

%%%%%%%%%%%%%%%%%%%%%%%%%

\subsubsection{A weak Laplacian and a weak gradient operator and its approximation}

Let $\mathcal{T}_{N}$ be a triangulation of the domain $\Omega$ as polygons in 2D, and $T$ denotes an element of the triangulation $\mathcal{T}_{N}$ and $\partial T$ its boundary.

We will define a weak Laplacian and a weak gradient operator on a space of weak functions.
Defining the space of local weak functions as:
\[
F(T) = \left\{u = \{u_{0},u_{b},\textbf{u}_{g}\} : u_{0} \in L^{2}(T), u_{b} \in H^{1/2}(\partial T), \textbf{u}_{g}\cdot \mathbf{n} \in H^{-1/2}(\partial T)\right\}.
\]
First, we define the set of normal directions on ${E}_N$ as follows:
\begin{align}\label{ch1_normaldirec}
{D}_h=\{ \textbf{n}_e: \; \textbf{n}_e \text{ is unit and normal to } e,\; e\in {E}_N\}
\end{align}
where $E_N$ denotes the set of all edges in $\mathcal{T}_{N}$.
Now given any non-negative integer $k$, we can define discrete weak function space on $T$ as:
\[
F_{N}(T) = \left\{u = \{u_{0},u_{b},u_{g}\textbf{n}_{e}\} : u_{0} \in \mathbb{P}_{k}(T), u_{b} \in \mathbb{P}_{k}(e), u_{g} \in \mathbb{P}_{k-1}(e),\; e\subset \partial T\right\},
\]
where $u_g$ can be interpreted as an approximation of  $\nabla u\cdot \boldsymbol{n}_e$.

Defining the weak finite element space globally, {\em i.e.}, on domain $\Omega$ as:
\[
S_N=S_{N}(\Omega) = \{u = \{u_{0},u_{b},u_{g}\textbf{n}_{e}\} : \{u_{0},u_{b},u_{g}\textbf{n}_{e}\}|\strut_T \in F_{N}(T), \forall \, T \in \mathcal{T}_{N}\}.
\]
Let $S_{N}^{0}(\Omega)$ be the subspace of $S_{N}$ such that:
\[
S_N^0=S_{N}^{0}(\Omega) = \{u = \{u_{0},u_{b},u_{g}\textbf{n}_{e}\} \in S_{N}(\Omega) : u_{b} = 0, u_{g}\vert_{e}= 0 \,\,\mbox{on}\, e \subset \partial T \cap \partial \Omega\}.
\]
Define $\mathcal{G}(T)$ as
\[
\mathcal{G}(T) = \{\vartheta : \vartheta \in H^{1}(T), \Delta \vartheta \in L^{2}(T)\}.
\]
Here, we can now think of a linear functional from the dual of $\mathcal{G}(T)$ that will actually be the weak local Laplacian that we denote as $\Delta_{w}$.

For any non-negative integer $m$,  the local weak Laplacian $\Delta_{w,m} u_N \in\mathbb{P}_m(T)$ of $u_N=\{ u_0,u_b,u_g\textbf{n}_e\}$ is defined as follows:
\begin{align}\label{ch1_weakLaplace}
(\Delta_{w,m}u,\psi)\strut_{T} = (u_{0},\Delta \psi)\strut_{T} -\langle u_{b}, \nabla \psi \cdot \mathbf{n}\rangle\strut_{\partial T} + \langle u_{g}\textbf{n}_e\cdot \mathbf{n}, \psi\rangle\strut_{\partial T}, \quad \forall \psi \in\mathbb{P}_m(T),
\end{align}
where $\langle \cdot , \cdot \rangle\strut_{\partial T}$ is the $L^{2}-$inner product on $\partial T$ and $\mathbf{n}$ represents the outward unit normal to $\partial T$.
Similarly, for any integer $l\geq 0$,  the weak gradient operator $\nabla_{w}u_N\in[\mathbb{P}_l(T)]^2$ can be defined on $T$ as:
\begin{align}\label{ch1_weakgradient}
(\nabla_{w}u_{N},\boldsymbol{\psi})\strut_{T} =-( u_0,\nabla \cdot \boldsymbol{\psi})_T+\langle u_b,\boldsymbol{\psi}\cdot \mathbf{n}\rangle_{\partial T} \quad \forall \boldsymbol{\psi} \in[\mathbb{P}_l(T)]^2.\end{align}

The discrete weak Laplacian operator on the WG finite element space $S_N$ is denoted by $\Delta_{w,k-2}$ and computed using (\ref{ch1_weakLaplace}) on each element $T$ for $k\geq 2$, that is,
\begin{align*}
\left.\left(\Delta_{w, k-2} u_N\right)\right|_T=\Delta_{w, k-2}\big(u_N|_T\big), \quad \forall u_N \in S_N.
\end{align*}
For the sake of simplicity, we do not write the subscript $k-2$ in $\Delta_{w, k-2}$.

For any $u_N=(u_0,u_b,u_g\textbf{n}_e)  \in S_{N}^{0}$ and $v_N=(v_0,v_b,v_g\textbf{n}_e)  \in S_{N}^{0}$  , we get our weak formulation as:
\begin{eqnarray}\label{ch1_bilinear}
\textsl{a}(u_N,v_N) = (\textsl{g}, v_{0}),
\end{eqnarray}
where
\begin{align}\label{ch1_bilinear_exp}
\textsl{a}(u_N,v_N) &= \sum_{T\in \mathcal{T}_N}(\varepsilon^2 \Delta_{w}u_N, \Delta_{w}v_N)\strut_T + \sum_{T\in \mathcal{T}_N}(\nabla_{w}u_N, \nabla_{w}v_N)\strut_T \\ \nonumber
&\quad+ \sum_{T\in \mathcal{T}_N}(a u_{0}, v_{0})\strut_T + s(u_N,v_N)
\end{align}
is the associated bilinear form, and
\begin{equation}\label{ch1_stab}
   s(u,v) = \sum_{T \in \mathcal{T}_{N}}^{}\left[ { \rho\strut_T}\langle\nabla u_{0}\cdot \textbf{n}_{e}-u_{g}, \nabla v_{0}\cdot \textbf{n}_{e}-v_{g}\rangle\strut_{\partial T}+ { \sigma\strut_T}\langle u_{0}-u_{b}, v_{0}-v_{b} \rangle\strut_{\partial T}\right]
\end{equation}
is the stabilizer term {with $\rho\strut_T$ and $\sigma\strut_T$ are  constants given by}
\begin{align}\label{ch1_sigma}
{ \rho_T}= \begin{cases} \varepsilon N, & \text { if } T \subset \Omega_0, \\[4pt]
\varepsilon N \ln ^{-1} N, & \text { if } T \subset \Omega \setminus \Omega_0 ,\end{cases}
\end{align}
 and
 \begin{align}\label{ch1_rho}
{ \sigma_T} = \begin{cases} N, & \text { if } T \subset \Omega_0, \\[4pt] 
\varepsilon^{-1}(N \ln ^{-1} N)^3, & \text { if } T \subset \Omega \setminus \Omega_0 .\end{cases}
 \end{align}

\begin{lemma}\label{triple_norm_lemma}
For any $u_N \in S^{0}_{N}$, let $\||u_N|\|$ be defined as
\[
\||u_N|\|^{2} = \textsl{a}(u_N,u_N).
\]
Then, $\|| \cdot |\|$ defines a norm in $S^{0}_{N}$.
\end{lemma}
\begin{proof}
{The proof can be found in \cite[Lemma 2.1]{cui2020uniform}.}
\end{proof}
%%%%%%%%%%%
\begin{lemma}
    The WG-FEM scheme \eqref{ch1_bilinear} admits at most one solution.
\end{lemma}
\begin{proof}
Let $u_N^1$ and $u_N^2$ be two solution of the WG-FEM approximation (\ref{ch1_bilinear}). Then, we have $w_N:= u_N^1-u_N^2\in S_N^0$ that satisfies
\begin{align}\label{ch1_unique-eq1}
    \textsl{a}(w_N,v_N)=0, \; \forall v_N=\{ v_0,v_b,v_g\textbf{n}_e\}\in S_N^0.
\end{align}
Taking $v_N=w_N$ in \eqref{ch1_unique-eq1} yields
\[ \|| w_N |\| = \textsl{a}(w_N,w_N)=0. \] The Lemma \ref{triple_norm_lemma} implies that $w_N=0$ leading to $u_N^1=u_N^2$. Thus, the proof is completed.
 \end{proof}
\begin{lemma}{(Trace-inverse estimate on triangle)}\label{trace}
Denote the edge of a mesh element $T$ by $e$, the length of the edge by $|e|$, and the area of the element by $|T|$. For any $\phi \in \mathbb{P}_{k}(T)$, for all elements $T \in \mathcal{T}_{N}$, the following trace inequality holds:
    \[
    \|\phi\|\strut_{e} \, \leq\,  \kappa\,\dfrac{|e|^{1/2}}{|T|^{1/2}}\|\phi\|\strut_{T},
    \]
    where $\kappa$ only depends on the degree of the polynomials $k$.
\end{lemma}
\begin{proof}
We can refer to \cite[Lemma 1.46]{di2011mathematical} for the proof.
\end{proof}

Introducing a discrete $H^{2}-$semi-norm on the space $S_{N}$ involving the usual derivatives as:
\begin{equation}\label{ch1_triplenorm}
    \||u|\|\strut_{M} = \left(\sum_{T\in \mathcal{T}_N}\Big(\varepsilon^2\Vert \Delta u_0\Vert\strut_T^2 + \Vert\nabla u_0\Vert\strut_T^2 + \Vert a u_{0}\Vert\strut_T^2\Big) + s(u_N,u_N)\right)^{\frac{1}{2}},
\end{equation}
where
\begin{align*}
    s(u_N,u_N) = \sum_{\partial T \in \mathcal{T}_{N}}^{}\rho\strut_T\|\nabla v_{0}-v_{g}\textbf{n}_{e}\|\strut_{\partial T}^2+ \sigma\strut_T\|v_{0}-v_{b}\|\strut_{\partial T}^2\,\,.
\end{align*}
 \begin{lemma}
{Let $\rho\strut_T$ and $\sigma\strut_T$ be given by \eqref{ch1_sigma} and \eqref{ch1_rho}, respectively. In Shishkin triangular mesh $\mathcal{T}_N$,}    there exist two positive constants $\kappa_1$ and $\kappa_2$ such that
\[      \kappa_1\||u|\|\strut_{M} \, \leq \, \||u|\| \, \leq \, \kappa_2\||u|\|\strut_{M}\,,      \] for all $u = \{u_{0},u_{b},{u}_{g}{\bf n_e}\} \in S_{N}$.
 \end{lemma}
%{\bf Proof.}
\begin{proof}
Using the definition of weak Laplacian \eqref{ch1_weakLaplace} for any $\{u_{0},u_{b},{u}_{g}{\bf n_e}\} \in S_{N}$ and integration by parts, we get {for $T\in \mathcal{T}_N$,}
\begin{align*}
(\Delta_{w}u,\psi)\strut_{T} &= (u_{0},\Delta \psi)\strut_{T} -\langle u_{b}, \nabla \psi \cdot \mathbf{n}\rangle\strut_{\partial T} + \langle u_{g}\textbf{n}_e\cdot \mathbf{n}, \psi\rangle\strut_{\partial T}, \quad \forall \psi \in\mathbb{P}_m(T),\\ \nonumber
&= -(\nabla u_{0}, \nabla \psi)\strut_{T} + \langle u_{0} - u_{b}, \nabla \psi \cdot \mathbf{n}\rangle\strut_{\partial T} + \langle u_{g}\textbf{n}_e\cdot \mathbf{n}, \psi \rangle\strut_{\partial T} \\ \nonumber
&= (\Delta u_{0}, \psi)\strut_{T} + \langle u_{0} - u_{b}, \nabla \psi \cdot \mathbf{n}\rangle\strut_{\partial T}+\langle (u_{g}\textbf{n}_e - \nabla u_{0})\cdot \mathbf{n}, \psi \rangle\strut_{\partial T}.
\end{align*}
Now, substituting $\psi = \Delta_{w} u$ into the above, we get
\begin{align}\label{ch1_weakdivnorm}
\|\Delta_{w} u\|\strut_{{T}}^2 = (\Delta u_{0}, \Delta_{w} u)\strut_{T} + \langle u_{0} - u_{b}, \nabla \left(\Delta_{w} u\right) \cdot \mathbf{n}\rangle\strut_{\partial T}+\langle (u_{g}\textbf{n}_e - \nabla u_{0})\cdot \mathbf{n}, \Delta_{w} u \rangle\strut_{\partial T} \cdot
\end{align}
{ It follows from Cauchy-Schwarz inequality and Lemma \ref{trace} that} %we get
\begin{eqnarray*}
    \varepsilon\|\Delta_{w} u\|\strut_{T}^{2} &\leq& \varepsilon\|\Delta u_{0}\|\strut_{T}\|\Delta_{w}u\|\strut_{T} + \varepsilon\sigma_T^{1/2}\|u_{0}-u_{b}\|\strut_{\partial T}\times\sigma_T^{-1/2}\|\nabla(\Delta_{w} u)\|\strut_{\partial T} \\[6pt] \nonumber
    & & + \varepsilon\rho_T^{1/2}\|(u_{g}\textbf{n}_e - \nabla u_{0})\cdot \mathbf{n}\|\strut_{\partial T}\times \rho_T^{-1/2}\|\Delta_{w} u\|\strut_{\partial T} \\ \nonumber
    & \leq& \varepsilon\|\Delta u_{0}\|\strut_{T}\|\Delta_{w}u\|\strut_{T} + \sigma_T^{1/2}\|u_{0}-u_{b}\|\strut_{\partial T} \times \\ \nonumber
    &&\, \varepsilon \sigma_T^{-1/2} \kappa\dfrac{|e|^{1/2}}{|T|^{1/2}}\|\nabla(\Delta_{w} u)\|\strut_{T} + \\ \nonumber
    &&  \rho_T^{1/2}\|(u_{g}\textbf{n}_e - \nabla u_{0})\cdot \mathbf{n}\|\strut_{\partial T}\times \varepsilon \rho_T^{-1/2}\kappa\dfrac{|e|^{1/2}}{|T|^{1/2}}\|\Delta_{w} u\|\strut_{T} \\ \nonumber
    & \leq & \varepsilon\|\Delta u_{0}\|\strut_{T}\|\Delta_{w} u\|\strut_{T} +\sigma_T^{1/2}\|u_{0}-u_{b}\|\strut_{\partial T} \times \\ \nonumber
    && \, \varepsilon \sigma_T^{-1/2} \kappa\dfrac{|e|^{1/2}}{|T|^{1/2}}h_{T}^{{-1}}\|\Delta_{w} u\|\strut_{T} + \\ \nonumber
    && \rho_T^{1/2}\|(u_{g}\textbf{n}_e - \nabla u_{0})\cdot \mathbf{n}\|\strut_{\partial T}\times \varepsilon \rho_T^{-1/2}\kappa\dfrac{|e|^{1/2}}{|T|^{1/2}}\|\Delta_{w} u\|\strut_{T},
\end{eqnarray*}
    which reduces to
    \begin{align}\label{ch1_divergenceEST}
    \varepsilon\|\Delta_{w} u\|\strut_{T}\, &\leq  \, C\left(\varepsilon\|\Delta u_{0}\|\strut_{T} + \sigma_T^{1/2}\|u_{0}-u_{b}\|\strut_{\partial T} + \rho_T^{1/2}\|(u_{g}\textbf{n}_e - \nabla u_{0})\cdot \mathbf{n}\|\strut_{\partial T}\right).
\end{align}
Here, we have the following inequalities: For $T\subset \Omega\setminus \Omega_0$,
\begin{align*}
  \frac{|e|^{1 / 2}}{|T|^{1 / 2}} &\leq C \varepsilon^{-1 / 2} N^{1 / 2} \ln ^{-1 / 2} N, \\
\rho_T^{-1 / 2} &\leq C \varepsilon^{-1/2} N^{-1 / 2} \ln ^{1 / 2} N,  
\end{align*}
and  for $T\subset  \Omega_0$,
\begin{align*}
  \frac{|e|^{1 / 2}}{|T|^{1 / 2}} &\leq C  N^{1 / 2},  \\
\rho_T^{-1 / 2} &\leq C \varepsilon^{-1/2} N^{-1 / 2} .  
\end{align*}
Thus, the coefficient $\varepsilon \rho_T^{-1/2}\kappa\cfrac{|e|^{1/2}}{|T|^{1/2}}$ does not include a negative power of $\varepsilon$. Furthermore, using the definition of $\sigma_T$ and the above estimates,  for any $T\subset \Omega_C$ with $\Omega_C:=\Omega_{41}\cup \Omega_{43}\cup\Omega_{21}\cup\Omega_{23}$, we have
\begin{align*}
{\varepsilon} \sigma_T^{-1/2} \kappa\dfrac{|e|^{1/2}}{|T|^{1/2}}h_{T}^{{-1}}&\leq \varepsilon\cdot \varepsilon^{1/2} (N^{-1}\ln N)^{3/2}\varepsilon^{-1 / 2} N^{1 / 2} \ln ^{-1 / 2} N \varepsilon^{-1 } (N^{-1} \ln  N)^{-1}  \\ &\leq C,  
\end{align*}
and for $T\subset \Omega\setminus(\Omega_0\cup\Omega_C)$, we have
\begin{align*}
{\varepsilon} \sigma_T^{-1/2} \kappa\dfrac{|e|^{1/2}}{|T|^{1/2}}h_{T}^{{-1}}&\leq \varepsilon\cdot \varepsilon^{1/2} (N^{-1}\ln N)^{3/2}\varepsilon^{-1 / 2} N^{1 / 2} \ln ^{-1 / 2} N \cdot N  \\ &\leq C,  
\end{align*}
where we have used that $\varepsilon N\leq 1$. Lastly, for $T\subset \Omega_0$, we have \begin{align*}
{\varepsilon} \sigma_T^{-1/2} \kappa\dfrac{|e|^{1/2}}{|T|^{1/2}}h_{T}^{{-1}}&\leq \varepsilon\cdot  N^{-1/2}\cdot N^{1 / 2} \cdot N  \\ &\leq C.  
\end{align*}
Similarly, using the definition of weak gradient \eqref{ch1_weakgradient} and integration by parts, we get
\begin{align}
(\nabla_{w}u,\boldsymbol{\psi})\strut_{T} &=-( u_0,\nabla \cdot \boldsymbol{\psi})\strut_T+\langle u_b,\boldsymbol{\psi}\cdot \mathbf{n}\rangle\strut_{\partial T} \quad \forall \boldsymbol{\psi} \in[\mathbb{P}_l(T)]^2\nonumber\\ \label{eq-1}
    & = (\nabla u_{0},\boldsymbol{\psi})\strut_{T} + \langle u_b - u_{0},\boldsymbol{\psi}\cdot \mathbf{n}\rangle\strut_{\partial T} \, .
\end{align}
Substituting $\boldsymbol{\psi} = \nabla_{w} u$ in the above equation, we have
\begin{align*}
    \|\nabla_{w}u\|\strut_{T}^{2} &= (\nabla u_{0},\nabla_{w}u)\strut_{T} + \langle u_b - u_{0},\nabla_{w}u\cdot \mathbf{n}\rangle\strut_{\partial T},
\end{align*}
which reduces to
\begin{align*}
     \|\nabla_{w}u\|\strut_{T}^{2} & \leq \|\nabla u_{0}\|\strut_{T}\|\nabla_{w}u\|\strut_{T} + \sigma\strut_T^{1/2}\|u_{0}-u_{b}\|\strut_{\partial T}\sigma\strut_T^{-1/2}\|\nabla_{w}u\|\strut_{\partial T}\\
     & \leq \|\nabla u_{0}\|\strut_{T}\|\nabla_{w}u\|\strut_{T} + \sigma\strut_T^{1/2}\|u_{0}-u_{b}\|\strut_{\partial T}\times\kappa\dfrac{|e|^{1/2}}{|T|^{1/2}}\sigma\strut_T^{-1/2}\|\nabla_{w}u\|\strut_{T},
\end{align*}
by using {Lemma \ref{trace}}.

We simplify the above inequality to obtain the following:
\begin{align}\label{ch1_gradientEST}
    \|\nabla_{w}u\|\strut_{T} \,\, \leq \,\, C \left(\|\nabla u_{0}\|\strut_{T} + \sigma\strut_T^{1/2}\|u_{0}-u_{b}\|\strut_{\partial T}\right) \cdot
\end{align}
Here, we have used the following observation: For any $T\subset \Omega\setminus\Omega_0$, one has
\begin{align*}
 \sigma_T^{-1/2} \kappa\dfrac{|e|^{1/2}}{|T|^{1/2}}&\leq  \varepsilon^{1/2} (N^{-1}\ln N)^{3/2}\varepsilon^{-1 / 2} N^{1 / 2} \ln ^{-1 / 2} N   \leq 1,  
\end{align*} and any $T\subset \Omega_0$, we have
\begin{align*}
 \sigma_T^{-1/2} \kappa\dfrac{|e|^{1/2}}{|T|^{1/2}}&\leq   N^{-1/2} N^{1 / 2}   \leq 1.  
\end{align*}
{In order to obtain} the lower bound, we choose $\boldsymbol{\psi} = \nabla u_{0}$ in \eqref{eq-1} to get
\begin{align*}
(\nabla_{w}u,\nabla u_{0})\strut_{T} = (\nabla u_{0},\nabla u_{0})\strut_{T} + \langle u_b - u_{0},\nabla u_{0}\cdot \mathbf{n}\rangle\strut_{\partial T}.
\end{align*}
{It follows from the Cauchy-Schwarz inequality and Lemma \ref{trace} that}
\begin{align*}
(\nabla u_{0},\nabla u_{0})\strut_{T}
& \leq \|\nabla_{w} u\|\strut_{T}\|\nabla u_{0}\|\strut_{T} + \sigma\strut_T^{1/2}\|u_0-u_b\|\strut_{\partial T}\times \sigma\strut_T^{-1/2}\|\nabla u_{0}\cdot\mathbf{n}\|\strut_{\partial T}\\
& \leq \|\nabla_{w} u\|\strut_{T}\|\nabla u_{0}\|\strut_{T} + \sigma\strut_T^{1/2}\|u_0-u_b\|\strut_{\partial T}\times \sigma\strut_T^{-1/2}\kappa\dfrac{|e|^{1/2}}{|T|^{1/2}}\|\nabla u_{0}\|\strut_{T},
\end{align*}
which reduces to
\begin{align}\label{ch1_gradv0bound}
    \|\nabla u_{0}\|\strut_{T} \,\, \leq \,\, C_{3}\left(\|\nabla_{w} u\|\strut_{T} + \sigma\strut_T^{1/2}\|u_0-u_b\|\strut_{\partial T}\right).
\end{align}
Using the bounds obtained in \eqref{ch1_divergenceEST} and \eqref{ch1_gradientEST}, we obtain
\begin{eqnarray*}
    \||u|\| &\leq& C_{1}\left(\sum_{T\in \mathcal{T}_{N}}^{}\varepsilon\|\Delta u_{0}\|\strut_{T} + \sigma\strut_T^{1/2}\|u_{0}-u_{b}\|\strut_{\partial T} + \rho\strut_T^{1/2}\|(u_{g}\textbf{n}_e - \nabla u_{0})\cdot \mathbf{n}\|\strut_{\partial T}\right)\\\nonumber
    && \,\,+ \,  C_{2}\left(\sum_{T\in \mathcal{T}_{N}}^{}\|\nabla u_{0}\|\strut_{T} + \sigma\strut_{T}^{1/2}\|u_{0}-u_{b}\|\strut_{\partial T}\right) + \sum_{T\in \mathcal{T}_{N}}^{}a_{2}\|u_{0}\|\strut_{T}^{2}\\
    & \leq & \kappa_{1} \||u|\|\strut_{M}.
\end{eqnarray*}
From \eqref{ch1_gradv0bound} and definition of $\||u|\|$, it is straightforward to have
\begin{align*}
     &\sum_{T \in \mathcal{T}_{N}}^{}\|\nabla u_{0}\|\strut_{T}^{2} \,\, \leq \,\, C\||u|\|^{2},\\
     &\sum_{T \in \mathcal{T}_N}^{} \sigma\strut_T\|u_0-u_b\|\strut_{\partial T}^{2} \,\, \leq \,\, C\||u|\|^2,\\
     &\sum_{T \in \mathcal{T}_N}^{} \rho\strut_T\|(\nabla u_0 - u_{g}\textbf{n}_e)\cdot \mathbf{n}\|\strut_{\partial T}^{2} \,\, \leq \,\, C\||u|\|^{2}.
\end{align*}
Now using integration by parts in \eqref{ch1_weakdivnorm}, we get
\begin{align*}
\|\Delta u_{0}\|\strut_{T}^{2} = (\Delta u_{0}, \Delta_{w} u)\strut_{T} - \langle u_{0} - u_{b}, \nabla \left(\Delta_{w} u\right) \cdot \mathbf{n}\rangle\strut_{\partial T}-\langle (u_{g}\textbf{n}_e - \nabla u_{0})\cdot \mathbf{n}, \Delta_{w} u \rangle\strut_{\partial T} \cdot
\end{align*}
Using trace inequality, inverse inequality and above bounds, we reach
\begin{eqnarray*}
\|\Delta u_{0}\|\strut_{T}^{2} &\leq& C\|\Delta u_{0}\|\strut_{T}\|\Delta_{w} u\|\strut_{T} \,\, ,  \\ [6pt]
\sum_{T \in \mathcal{T}_N}\|\Delta u_{0}\|\strut_{T} & \leq & C\|\Delta_{w} u\|\strut_{T}\,\,.
\end{eqnarray*}
which will give the lower bound for $\||u|\|$, {\em i.e.},
\[
\kappa_{2}\||u|\|\strut_{M} \,\, \leq  \,\,\||u|\| \cdot
\]
This completes the proof.
\end{proof}
%%%%%%%%%%%%%%%%%%%%%%%%%%%%%%%%%%%%%%%%%%%%%%%%%%%%%%%%%%%%%%%%%%%%%%%%%%%%%%%%%%%

\section{Local discretization of weak Galerkin FEM}\label{WG_Discre}

\subsection{Computation of local stiffness matrices}

Let $N_{0}, N_{b} \,\, \mbox{and} \,\, N_{g}$ represent the dimensions for the polynomial spaces of $u_{0}, u_{b}, u_{g}$, respectively, and let $\{\psi_{0,i}\}\strut_{i=1}^{N_{0}} \,\,, \{\psi_{b,i}\}\strut_{i=1}^{N_{b}} \,\,, \{\psi_{g,i}\}\strut_{i=1}^{N_{g}}$ be their corresponding bases. Then $u$ can be expressed by
\[
u|\strut_{T} = \left\{\displaystyle\sum_{i=1}^{N_{0}} u_{0,i}\psi_{0,i}, \sum_{i=1}^{N_{b}}u_{b,i}\psi_{b,i},\sum_{i=1}^{N_{g}}u_{g,i}\psi_{g,i}\right\}.
\]
On each element $T$, the local stiffness matrix $S_{T}$ for \eqref{ch1_bilinear} can then be expressed as
\begin{eqnarray}\label{ch1_local_stiffness}
S_{T} = \begin{bmatrix}
S_{0,0} & S_{0,b} & S_{0,g}\\
S_{b,0} & S_{b,b} & S_{b,g}\\
S_{g,0} & S_{g,b} & S_{g,g}
\end{bmatrix}
\end{eqnarray}
These matrices are defined as:
\[
S_{0,0} = [\textsl{a}(\psi_{0,i}, \psi_{0,j})]_{i,j}, \quad S_{0,b} = [\textsl{a}(\psi_{0,i}, \psi_{b,j})]_{i,j}, \quad S_{0,g} = [\textsl{a}(\psi_{0,i}, \psi_{g,j})]_{i,j}
\]
\[
S_{b,0} = [\textsl{a}(\psi_{b,i}, \psi_{0,j})]_{i,j}, \quad S_{b,b} = [\textsl{a}(\psi_{b,i}, \psi_{b,j})]_{i,j}, \quad S_{b,g} = [\textsl{a}(\psi_{b,i}, \psi_{g,j})]_{i,j}
\]
\[
S_{g,0} = [\textsl{a}(\psi_{g,i}, \psi_{0,j})]_{i,j}, \quad S_{g,b} = [\textsl{a}(\psi_{g,i}, \psi_{b,j})]_{i,j}, \quad S_{g,g} = [\textsl{a}(\psi_{g,i}, \psi_{g,j})]_{i,j}
\]
where $i,j$ correspond to the row and column numbers and $\textsl{a}(\cdot ,\cdot)$ is defined as earlier in \eqref{ch1_bilinear_exp}.

Before we can calculate $S_{T}$, we need to evaluate the discrete Laplacian operator $\Delta_{w}$ and discrete gradient operator $\nabla_{w}$.

Let $\psi_{w,i}, i = 1,\ldots , N_{w}$ be a set of basis functions of $P_{r}(T) \subset \mathcal{G}(T)$. Then we can express $\Delta_{w}v$ as follows:
\[
\Delta_{w}v|\strut_{T} = \sum_{i=1}^{N_{w}} (\Delta_{w}v)_{i} \psi_{w,i}.
\]
Then $\Delta_{w}v$ corresponding to the bases functions $v_{i}$$=$$\{\psi_{0,i}\}\strut_{i=1}^{N_{0}}, \{\psi_{b,i}\}\strut_{i=1}^{N_{b}}, \{\psi_{g,i}\}\strut_{i=1}^{N_{g}}$ can be expressed as
\[
\begin{bmatrix}
\Delta_{w}\psi_{0,i}\\
\Delta_{w}\psi_{b,i}\\
\Delta_{w}\psi_{g,i}
\end{bmatrix}
=
\begin{bmatrix}
C_{0}\\
C_{b}\\
C_{g}
\end{bmatrix}
[\psi_{w,i}]
\]
where $C_{0}, C_{b}, C_{g}$ are matrices of coefficients whose values are to be evaluated. More details can be found in \cite{burkardt2020high}.

Now, using the definition of the discrete weak gradient, for $u|\strut_{T}$ we can find the vector form of $\nabla_{w}v$ on $T$ by
\begin{equation}\label{ch1_vgrad}
\mathcal{D}|\strut_{T}(\nabla_{w} v) = -\mathcal{E}_{T}\psi_{0} + \mathcal{F}_{T}\psi_{b} ,
\end{equation}
where $\mathcal{D}_{T}$ is $N_{g}\times N_{g}$ matrix, $\mathcal{E}_{T}$ is $N_{g}\times N_{0}$ matrix and $\mathcal{F}_{T}$ is $N_{g}\times N_{b}$ matrix. Exact representation of these matrices can be found in \cite{mu2013computational}. Once these matrices are evaluated, by using \eqref{ch1_vgrad} we can compute the weak gradient of bases functions $ \psi_{0,i}$ and $\psi_{b,i}$ on $T$ as follows:
\[
(\nabla_{w}\psi_{0,i}) = -\mathcal{D}^{-1}_{T}\mathcal{E}_{T}{\bf e}^{N_{0}}_{i}, \quad (\nabla_{w}\psi_{b,i}) = \mathcal{D}^{-1}_{T}\mathcal{F}_{T}{\bf e}^{N_{b}}_{i},
\]
where ${\bf e}^{N_{0}}_{i}$ and ${\bf e}^{N_{b}}_{i}$ are vectors of order $N_{0} \times 1$ and $N_{b} \times 1$ respectively, with the property that their $i$th entry is unity and others zero.

Define matrices
\begin{eqnarray*}
A_{T} &=& \left[(\varepsilon^2\psi_{w,i},\psi_{w,j})\strut_{T}\right]_{i,j} ,\\
B_{T} &=& \left[(\psi_{g,i},\psi_{g,j})\strut_{T}\right]_{i,j} ,\\
C_{T} &=& \left[(a\psi_{0,i},\psi_{0,j})\strut_{T}\right]_{i,j} ,
\end{eqnarray*}
where $(\cdot , \cdot)_{T}$ denote the standard $L^{2}(T)$ or $[L^{2}(T)]^{d}$ inner product, whichever is appropriate.

Now, we are in a position, where we can find the local stiffness matrix $S_{T}$ and the same is specified in the following lemma:

\begin{lemma}
We can compute $S_{T}$ (the local stiffness matrix) as defined in \eqref{ch1_local_stiffness} by using the following formula:
\begin{eqnarray*}
S_{0,0} &=& C_{0}A_{T}C^{t}_{0} + \mathcal{E}^{t}_{T}\mathcal{D}^{-t}_{T}B_{T}\mathcal{D}^{-1}_{T}\mathcal{E}_{T} + C_{T},\\
S_{0,b} &=& C_{0}A_{T}C^{t}_{b} - \mathcal{E}^{t}_{T}\mathcal{D}^{-t}_{T}B_{T}\mathcal{D}^{-1}_{T}\mathcal{F}_{T},\\
S_{0,g} &=& C_{0}A_{T}C^{t}_{g},\\
S_{b,0} &=& C_{b}A_{T}C^{t}_{0} - \mathcal{F}^{t}_{T}\mathcal{D}^{-t}_{T}B_{T}\mathcal{D}^{-1}_{T}\mathcal{E}_{T},\\
S_{b,b} &=& C_{b}A_{T}C^{t}_{b} + \mathcal{F}^{t}_{T}\mathcal{D}^{-t}_{T}B_{T}\mathcal{D}^{-1}_{T}\mathcal{F}_{T},\\
S_{b,g} &=& C_{b}A_{T}C^{t}_{g},\\
S_{g,0} &=& C_{g}A_{T}C^{t}_{0},\\
S_{g,b} &=& C_{g}A_{T}C^{t}_{b},\\
S_{g,g} &=& C_{g}A_{T}C^{t}_{g},
\end{eqnarray*}
where $t$ represents the standard matrix transpose.
\end{lemma}
%%%%%%%%%%%%%%%%%%%%%%%%%%%%%%%%%%%%%%%%%%%%%%%%%%%%%%%%%%%%%%%%%%%%%%%%%%%%%%%%%%%%%
\subsection{Local discrete stabilizer term}
The local stabilizer term which we have taken is
\begin{eqnarray*}
s(u,v) := \sum_{T \in \mathcal{T}_N}^{}\big[\rho\strut_T\langle\nabla u_{0}-u_{g}\textbf{n}_{e}, \nabla v_{0}-v_{g}\textbf{n}_{e}\rangle\strut_{\partial T} + \sigma\strut_T\langle u_{0}-u_{b}, v_{0}-v_{b} \rangle\strut_{\partial T}\big].
\end{eqnarray*}
Now let
\begin{eqnarray*}
s_{1} &=& \rho\strut_T\langle\nabla u_{0}-u_{g}\textbf{n}_{e}, \nabla v_{0}-v_{g}\textbf{n}_{e}\rangle\strut_{\partial T}, \\ [6pt]
s_{2} &= & \sigma\strut_T\langle u_{0}-u_{b}, v_{0}-v_{b} \rangle\strut_{\partial T}, \end{eqnarray*}
then, we can rewrite $s(u,v)$ as
\begin{equation*}
    s(u,v) = \sum_{T \in \mathcal{T}_N} (s_1 + s_2).
\end{equation*}
Now, by using the form of $u|\strut_{T}$ as defined in the beginning of this section and substituting it in $s_1$ and $s_2$, we get the matrix corresponding to the discrete stabilizer term on element $T$.
	
\section{Error Analysis }\label{Error_Analysis}

In this section, we will present error estimates for the Weak Galerkin finite element solution $u_{N}$ of \eqref{ch1_bilinear}.
\subsection{Interpolation operator }
In this subsection, four interpolation operators $I_0, I_b, {I}_g$ and $I_h$  are introduced for the interpolation error analysis. For each $T\in \mathcal{T}_N$,
 $I_0:H^2(T)\to \mathbb{P}_k(T)$ is the Lagrange interpolation, $I_b:L^2(e)\to \mathbb{P}_k(e)$ is the $L^2$-projection onto $\mathbb{P}_k(e)$, and ${I}_g$  is the $L^2$-projection onto $\mathbb{P}_{k-1}(e)$. Finally, for any $u\in H^2(\Omega)$, we define the projection $I_h$ into the finite element space $S_N$ such that $I_h u=\{I_0 u,I_b u, {I}_g (\nabla u\cdot \mathbf{n}_e)\mathbf{n}_e\}$ on the element $T$.
 \begin{remark}
 Since the Shishkin mesh is highly anisotropic, we need sharp anisotropic interpolation error estimates inside of the element. For this reason, we introduce the Lagrange interpolation $I_0$ instead of the $L^2$-projection. The $L^2$-projections $I_b$ and ${I}_g$ are the commonly used projections in the error analysis of WG-FEM for the biharmonic problems, see \cite{cui2020uniform}, \cite{zhang2015weak}.
 \end{remark}

 We need the following local anisotropic interpolation error estimate for our analysis. But before that let us state the maximal angle condition \cite{local}.
 
{\bf Maximal angle condition}: There is a constant $\gamma_* \leq \pi$ (independent of $h_T$ and $T\in\mathcal{T}_N$) such that the maximal interior angle $\gamma\strut_T$ of any element $T$ is bounded by $\gamma\strut_* $.

\begin{theorem}[\cite{local}]\label{anistropic}
    Assume that the element $T$ fulfills the maximal angle condition. Taking $\gamma=\left(\gamma_1, \gamma_2\right)$ as a multiindex where $m:=\gamma_1+\gamma_2$ and $u \in \mbox{\escript{C}}(\bar{T})$ are functions for which $D^\gamma u \in H^{\bar\ell-m}(T)$ where $\bar\ell, m \in \mathbb{N}$ such that $0 \leq m <\bar\ell$ and $2 \leq \bar\ell \leq k+1$ hold. Then the following estimate:
\begin{equation*}
\left\|D^\gamma\left(u-I_0 u\right)\right\|_{L^2(T)} \leq C \sum_{|\beta|=1} h^{(\bar\ell-m) \beta}\left\|D^{\gamma+(\bar\ell-m) \beta} u\right\|_{L^2(T)}
\end{equation*}
holds with $h^\beta:=h_{x,T}^{\beta_1} h_{y,T}^{\beta_2}$.
 \end{theorem}
\begin{proof}
    For the detailed proof, we refer the reader to \cite[Theorem 2.1]{local}.
\end{proof}

In our analysis, the following interpolation error estimates on the edge $e\subset \partial T$ for each $T\in \mathcal{T}_N$ will be helpful.
\begin{lemma}\label{edgeinterpolation}
For any edge $e\subset \partial T$, the interpolation $I_0$ and the $L^2$-projections $I_b$ and ${I}_g$ satisfy
    \begin{eqnarray*}
\left\|{I}_0 u-u\right\|_{L^{\infty}(e)}  &\leq& \left\{\begin{array}{rll}
C N^{-s}+C \varepsilon N^{-\lambda},  &\text { if } e \subset \Omega_0, \\[6pt]
C \varepsilon^{3/2} (N^{-1} \ln  N)^{s},  &\text { if } e \subset \Omega\setminus \Omega_0,
\end{array}\right. \\[6pt]
\left\|\nabla\left({I}_0 u-u\right)\right\|_{L^{\infty}(e)}&\leq&\left\{\begin{array}{rll}
C N^{-(s-1)}+C  N^{-\lambda}, & \text { if } e \subset \Omega_0, \\[6pt]
C \varepsilon^{1/2} N^{-(s-1)} \ln ^{s-1} N, & \text { if } e \subset \Omega\setminus \Omega_0,
\end{array}\right. \\[6pt]
\varepsilon \left\|\Delta\left({I}_0 u-u\right)\right\|_{L^{\infty}(e)}& \leq& \left\{\begin{array}{rll}
C\min\{ \varepsilon^{1/2} N^{2}{,} \,N^{1/2}\} N^{-\lambda}, & \text { if } e \subset \Omega_0, \\[6pt]
C \varepsilon^{1/2} N^{-(s-2)} \ln ^{s-2} N, &  \text { if } e \subset \Omega\setminus \Omega_0,
\end{array}\right. \\[6pt]
\varepsilon^2 \left\|\nabla( \Delta\left({I}_0 u-u\right))\right\|_{L^{\infty}(e)}& \leq& \left\{\begin{array}{rll}
C\min\{ \varepsilon^{1/2} N^{3}{,}\, N\} N^{-\lambda}, & \text { if } e \subset \Omega_0, \\[6pt]
C \varepsilon^{1/2} N^{-(s-3)} \ln ^{s-3} N, &  \text { if } e \subset \Omega\setminus \Omega_0,
\end{array}\right. \\[6pt]
\left\|{I}_b u-u\right\|_{L^{\infty}(e)}& \leq& \left\{\begin{array}{rll}
C N^{-s}+C \varepsilon N^{-\lambda}, & \text { if } e \subset \Omega_0, \\[6pt]
C \varepsilon^{3/2} (N^{-1} \ln  N)^{s}, & \text { if } e \subset \Omega\setminus \Omega_0,\end{array}\right. \\[6pt]
\left\|({I}_g (\nabla u\cdot \mathbf{n}_e)\mathbf{n}_e-\nabla u\right\|_{L^{\infty}(e)} & \leq& \left\{\begin{array}{rll}
C N^{-(s-1)}+C  N^{-\lambda}, & \text { if } e \subset \Omega_0, \\[6pt]
C \varepsilon^{1/2} N^{-(s-1)} \ln ^{s-1} N, & \text { if } e \subset \Omega\setminus \Omega_0,
\end{array}\right.
\end{eqnarray*}
where $u\in H^s(\Omega)$ for $3=s=k+1$ or  $4\leq s\leq k+1$.
\end{lemma}
\begin{proof}
    For any $e\subset\partial T$, $T\in \mathcal{T}_N$,  it follows from the fact $u, I_0 u\in \mbox{\escript{C}}^{3}(\overline{T})$ that
    \begin{eqnarray*}
         \left\|{I}_0 u-u\right\|_{L^{\infty}(e)}&\leq& \left\|{I}_0 u-u\right\|_{L^{\infty}(T)},\\
         \left\|\nabla({I}_0 u-u)\right\|_{L^{\infty}(e)}&\leq& \left\|\nabla({I}_0 u-u)\right\|_{L^{\infty}(T)},\\
          \left\|\Delta({I}_0 u-u)\right\|_{L^{\infty}(e)}&\leq& \left\|\Delta({I}_0 u-u)\right\|_{L^{\infty}(T)}.
    \end{eqnarray*}

Now,  from Theorem \ref{anistropic} and Lemma \ref{solutiondecom}, we have the error estimates for $\left\|{I}_0 u-u\right\|_{L^{\infty}(T)}$,  $\left\|\nabla ({I}_0 u-u)\right\|_{L^{\infty}(T)}$  and  $\left\|\Delta({I}_0 u-u)\right\|_{L^{\infty}(T)}$ which give the desired estimates for $\left\|{I}_0 u-u\right\|_{L^{\infty}(e)}$,  $\left\|\nabla({I}_0 u-u)\right\|_{L^{\infty}(e)}$ and $\left\|\Delta({I}_0 u-u)\right\|_{L^{\infty}(e)}$. A similar argument was given in \cite{zhanguniform}.

Now we follow the ideas used in \cite{franz2014c0} to derive the bound for $\varepsilon^2 \left\|\nabla( \Delta\left({I}_0 u-u\right))\right\|_{L^{\infty}(T)}$ for any $ T\in \mathcal{T}_N$.
The desired estimate follows from the following fact:
\[  \varepsilon^2 \left\|\nabla( \Delta\left({I}_0 u-u\right))\right\|_{L^{\infty}(e)}\leq \varepsilon^2 \left\|\nabla( \Delta\left({I}_0 u-u\right))\right\|_{L^{\infty}(T)}.\]
Using the decomposition of solution and the linearity of the interpolation operator $I_0$, we split the interpolation error into \[ u-I_0u= \mathscr{S}-I_0\mathscr{S} + \sum_{i \in K}^{}(\mathscr{E}_{i}-I_0\mathscr{E}_i). \]
Let $\eta:=\mathscr{S}-I_0\mathscr{S}$.  Since the bounds on $\mathscr{S}$ are independent of $\varepsilon$, we invoke Theorem \ref{anistropic} for $m =3$ and $m<s$ for all elements $T \in \mathcal{T}_N$ with $h_{x,T}, \; h_{y,T} \leq h_2$ and obtain
\begin{eqnarray}\label{smooth}
\Vert \eta^{(m)}\Vert\strut_{L^\infty(T)}^2 & \leq & C \sum_{|\beta|=m}\Bigg(h_x^{s-m}\left\|D^\beta \frac{\partial^{s-m} \mathscr{S}}{\partial x^{s-m}}\right\|_{L^\infty(T)} +h_y^{s-m}\left\|D^\beta \frac{\partial^{s-m} \mathscr{S}}{\partial y^{s-m}}\right\|_{L^\infty(T)}\Bigg)^2 \nonumber \\ [2pt]
& \leq & C h_2^{2(s-m)}\Vert \mathscr{S}^{(s)}\Vert\strut_{L^\infty(T)}^2,
\end{eqnarray}
and thus, we have
\begin{equation*}
\sum_{T\in \mathcal{T}_N} \Vert \eta^{(m)}\Vert\strut_{L^\infty(T)}^2 \leq  C N^{2(m-s)} .
\end{equation*}

When $s=m=3$, the polynomial degree $k=2$ implies $|\eta|\strut_{H^3(T)}=|\mathscr{S}|\strut_{H^3(T)}$ and we have the above estimates as well.

Next, define  $\zeta:=\mathscr{E}_1-I_0\mathscr{E}_1$. Knowing that  on the region $\Omega_1^*=\Omega_{41} \cup \Omega_4 \cup \Omega_{43}$ the  boundary layer $\mathscr{E}_1$ is resolved, the rest of the proof will consist of  two cases. The interpolation errors of the other layers can be proved in a similar way.

 Let $T \subseteq \Omega_1^*$.  We shall estimate  $\Vert \zeta^{(3)}\Vert\strut_{L^\infty(T)}$. Using Theorem \ref{anistropic} with  $m=s=3=k+1$, we have
$$
\begin{aligned}
\Vert \zeta^{(3)}\Vert\strut_{L^\infty(T)}^2 & \leq C \sum_{|\beta|=m}\Bigg(h_{x,T}^{2(s-m)}\left\|D^\beta \frac{\partial^{s-m} \mathscr{E}_1}{\partial x^{s-m}}\right\|_{L^\infty(T)}^2\\
&\quad\quad\quad\quad\quad\quad+h_{y,T}^{2(s-m)}\left\|D^\beta \frac{\partial^{s-m} \mathscr{E}_1}{\partial y^{s-m}}\right\|_{L^\infty(T)}^2\Bigg) \\
& =C\left(h_{x,T}^{2(s-m)}\left\Vert\frac{\partial^{s-m} \mathscr{E}_1^{(m)}}{\partial x^{s-m}}\right\Vert_{L^\infty(T)}^2+h_{y,T}^{2(s-m)}\left\Vert\frac{\partial^{s-m} \mathscr{E}_1^{(m)}}{\partial y^{s-m}}\right\Vert _{L^\infty(T)}^2\right) .
\end{aligned}
$$
Noting that $h\strut_{x,T}=h_1$ on $\Omega_1^* \backslash \Omega_4$ and summing up the  estimate above gives
\begin{equation*}
\sum_{T \in \mathcal{T}_N\left(\Omega_1^*\right)}\Vert\zeta^{(3)}\Vert\strut_{L^\infty(T)}^2 \leq C\left(h_1^{2(s-m)}\left\Vert \mathscr{E}_1^{(m)}\right\Vert_{L^\infty\left(\Omega_1^*\right)}^2+h_2^{2(s-m)}\left\Vert \frac{\partial^{s-m} \mathscr{E}_1^{(m)}}{\partial y^{s-m}}\right\Vert_{L^\infty\left(\Omega_1\right)}^2\right) \text {. }
\end{equation*}
 From \eqref{ch1_Dsol3} and \eqref{ch1_Dsol4} and using the fact $\varepsilon h_1^{-1}\leq C h_2^{-1}$ with the choice of $\lambda$, we have for $m =3$:
\begin{align}\label{secondone}
\begin{aligned}
 \sum_{T \subset \Omega_1^* }\Vert \zeta^{(m)}\Vert\strut_{L^\infty(T)}^2&=\sum_{T \subset\Omega_1^*}\Vert\zeta^{(m)}\Vert\strut_{L^\infty(T)}^2 \leq C \varepsilon^{3-2 m}\left(h_1^{2(s-m)} \varepsilon^{2 m-2 s}+h_2^{2(s-m)}\right) \\
& \leq C \varepsilon^{3-2 m}\left(h_1 \varepsilon^{-1}\right)^{2(s-m)} \\
& =C \varepsilon^{3-2 m}\left(N^{-1} \ln N\right)^{2(s-m)} \text {. } \\
&
\end{aligned}
\end{align}
On $\Omega_1^c:=\Omega \backslash \Omega_1^*$, the layer $\mathscr{E}_1$ is not resolved. Imitating the techniques used in \cite{franz2014c0}, we can prove the following estimates:
\begin{eqnarray}
    \label{estimate-1}
    \sum_{T \subset \Omega_1^c }\Vert \zeta^{(m)}\Vert\strut_{L^\infty(T)}^2&\leq  C \varepsilon^{3-2 m}\left(h_2^{-m} \exp \left(-\dfrac{\lambda}{\varepsilon}\right)\right)^2,\\\label{estimate-2}
    \sum_{T \subset \Omega_1^c }\Vert \zeta^{(m)}\Vert\strut_{L^\infty(T)}^2&\leq  C \varepsilon^{1-2 m}\left(h_2 \exp \left(-\dfrac{\lambda}{\varepsilon}\right)\right)^2,\\\label{estimate-3}
    \sum_{T \subset \Omega_1^c }\Vert \zeta^{(m)}\Vert\strut_{L^\infty(T)}^2&\leq  C\varepsilon^{2-2 m} h_2^{1-m}\left(\exp \left(-\dfrac{\lambda}{\varepsilon}\right)\right)^2.
\end{eqnarray}
Using the estimates \eqref{estimate-1}--\eqref{estimate-3}, we conclude that
\begin{align}\label{lastestimate}
    \sum_{T \subset \Omega_1^c }\Vert \zeta^{(3)}\Vert\strut_{L^\infty(T)}^2&\leq  C\min\{ \varepsilon^{-3}h_2^{-6},\varepsilon^{-4}h_2^{-2},\varepsilon^{-5}h_2^2\}\left(\exp \left(-\frac{\lambda}{\varepsilon}\right)\right)^2.
\end{align}
Consequently, choosing one of these minima yields the following desired estimate:
\begin{align*}
   \sum_{T\in \mathcal{T}_N} \varepsilon^4 \Vert \nabla( \Delta\left({I}_0 u-u\right))\Vert\strut_{L^{\infty}(T)}^2&=\varepsilon^4 \Vert\nabla( \Delta\left({I}_0 u-u\right))\Vert\strut_{L^{\infty}(\Omega_1^*)}^2\\
   &\qquad+\varepsilon^4 \Vert\nabla( \Delta\left({I}_0 u-u\right))\Vert\strut_{L^{\infty}(\Omega\setminus \Omega_1^*)}^2\\
   &\leq C\Big( \varepsilon ^{1/2}\left(N^{-1} \ln N\right)^{(s-3)}\\
   &\qquad+\min\{ \varepsilon^{1/2} N^{3}{,} N\} N^{-\lambda}\Big)^2.
\end{align*}

The $L^\infty-$stability of the $L^2-$projections on the Shishkin mesh \cite{oswald2012} give that
\begin{align*}
\left\|{I}_b u-u\right\|_{L^{\infty}(e)} \leq C\left(\left\|u^I-u\right\|_{L^{\infty}(e)}+ \left\|{I}_b (u^I-u)\right\|_{L^{\infty}(e)}\right)\leq C \left\|u^I-u\right\|_{L^{\infty}(e)}
\end{align*}
and
\begin{align*}
\left\|({I}_g (\nabla u\cdot \mathbf{n}_e)\mathbf{n}_e-\nabla u\right\|_{L^{\infty}(e)} &\leq C\Big(\left\|\nabla(u^I-u)\right\|_{L^{\infty}(e)}\\
&\qquad+ \left\|({I}_g (\nabla (u^I-u))\cdot \mathbf{n}_e)\mathbf{n}_e\right\|_{L^{\infty}(e)}\Big),\\[4pt] 
&\leq C \left\| \nabla(u^I -u)\right\|_{L^{\infty}(e)},
\end{align*}
where $u^I$ is the Lagrange interpolation of $u$ on  the edge $e\subset \partial T$ for a given $T\in \mathcal{T}_N$. Since $u, I_0 u\in \mbox{\escript{C}}^{1}(\overline{T})$, we have $u^I\vert_e=I_0 u\vert_e$ and $\nabla u^I\vert_e= \nabla(I_0 u)\vert_e$  for any $e\subset \partial T$.  Thus, the last two error estimates follow from the first two bounds, which are the required results.
\end{proof}

The following lemma from \cite[Lemma 3.5]{franz2014c0} will be used in the error analysis:
\begin{lemma}\cite{franz2014c0}\label{interpolationlemma}
Let the solution of the problem \eqref{ch1_md_problem} $u$ belongs to $H^s(\Omega)$ for
$s=3=k+1$ or $4 \leq s \leq k+1$ and satisfies Lemma \ref{solutiondecom}. Then on a family of Shishkin meshes as defined by \eqref{ch1_h1_h2}-\eqref{ch1_tensor_mesh}, the following estimates hold:
\begin{eqnarray*}
 \varepsilon^2 \sum_{T \in \mathcal{T}_N}|u-I_0 u|\strut_{H^2(T)}^2 &\leq& C\left(\varepsilon^{\frac{1}{2}}\left(N^{-1} \ln N\right)^{s-2}+\min \left\{\varepsilon^{\frac{1}{2}} N^2, N^{\frac{1}{2}}\right\} N^{-\lambda}\right)^2, \\
|u-I_0 u|\strut_{H^1(\Omega)}^2 &\leq& C\left(\varepsilon^{\frac{1}{2}}\left(N^{-1} \ln N\right)^{s-1}+N^{-(s-1)}+N^{-\lambda}\right)^2, \\
 \|u-I_0 u\|\strut_{L^2(\Omega)}^2 &\leq& C\left(\varepsilon^{\frac{3}{2}}\left(N^{-1} \ln N\right)^s+N^{-s}+\varepsilon N^{-\lambda}\right)^2 .
\end{eqnarray*}
The constants $C$ are independent of $\varepsilon$ and mesh properties.
\end{lemma}

\subsection{Error equations}
The proposed WG-FEM has a consistency error since the exact solution does not satisfy the numerical scheme.  This consistency error deteriorates the classical orthogonal property of the classical FEM and some discontinuous Galerkin methods.  In order to estimate this consistency error, we will derive some error equations.

\begin{lemma}\label{errorequation}
Let $I_h u$  be the interpolation of the solution  $u$ of \eqref{ch1_md_problem}. Then we have the following identities:
\begin{eqnarray}\label{ch1_2Laplace}
   \varepsilon^2(\Delta^2 u,v_0) &=& \varepsilon^2(\Delta_w(I_h u),\Delta_w v_N)-Z_1(u,v_N), \quad \forall v_N \in S^{0}_{N},\\[2pt] \label{ch1_Laplace}
-\left(\Delta u, v_0\right) &=& \left(\nabla_w\left(I_h u\right), \nabla_w v_N\right)-Z_2\left(u, v_N\right),\\[2pt] \label{ch1_reaction}
 \left(au, v_0\right) &=& \left(a I_h u, v_0\right)-Z_3(u,v_N),
\end{eqnarray}
where
\begin{eqnarray} \label{ch1_Z1}
Z_{1}(u,v_N) &=& \varepsilon^2\big( \Delta (-u+I_0 u),\Delta v_0\big)+\varepsilon^2\big\langle \Delta u-\Delta I_0 u, (\nabla v_0-v_g\textbf{n}_e) \cdot \textbf{n}\big\rangle \nonumber\\[2pt]
&& \, -\, \varepsilon^2\big\langle \nabla(\Delta u-\Delta I_0 u)\cdot \textbf{n},v_0-v_b\big\rangle \nonumber \\[2pt] 
&& \, +\, \varepsilon^2\langle ( {I}_{g}(\nabla u)\mathbf{n_e}-\nabla(I_0 u))\cdot  {n},\Delta _w v_N\rangle \nonumber \\[2pt]
&& \, -\, \varepsilon^2\langle I_bu-I_0 u, \nabla (\Delta _w v_N)\cdot  {n} \rangle , \\[2pt] \label{ch1_Z2}
Z_{2}(u,v_N) &=& \left(\nabla\left(u-I_0 u\right), \nabla v_0\right)+ \left\langle\left(\nabla u-\nabla I_0 u\right) \cdot \mathbf{n}, v_b-v_0\right\rangle \nonumber\\[2pt]
&& \, +\,\big\langle I_0 u-I_b u, \nabla_w v_N\cdot \mathbf{n}\big\rangle ,\\[2pt] \label{ch1_Z31}
Z_{3}(u,v_N) &=& (a(I_0 u - u),v_0).
\end{eqnarray}
\end{lemma}
\begin{proof}
For $v_N=\{v_0,v_b,v_g\textbf{n}_e\}\in S_N^0$, integration by parts gives
\begin{align}
    \label{ch1_z1-1}
    \begin{split}
        \big( \Delta (I_0 u),\Delta v_0\big)&= \big( \Delta (I_0 u-u),\Delta v_0\big)+\big( \Delta u, \Delta v_0\big),\\[2pt]
        &=\big( \Delta (I_0 u-u),\Delta v_0\big)+\big( \Delta^2 u, v_0\big)+\big\langle \Delta u, \nabla v_0 \cdot \textbf{n}\big\rangle\\[2pt]
        &\qquad-\big\langle \nabla(\Delta u)\cdot \textbf{n},v_0\big\rangle,\\[4pt]
        &= \big( \Delta (I_0 u-u),\Delta v_0\big)+\big( \Delta^2 u, v_0\big)+\big\langle \Delta u, (\nabla v_0-v_g\textbf{n}_e) \cdot \textbf{n}\big\rangle\\[2pt]
        &\qquad-\big\langle \nabla(\Delta u)\cdot \textbf{n},v_0-v_b\big\rangle,
    \end{split}
\end{align}
where we have used the fact that $v_g$ and $v_b$ vanish on the boundary of the domain.

It follows from the definition of the weak Laplacian operator (\ref{ch1_weakLaplace}), integration by parts  and (\ref{ch1_z1-1}) that for any $v_N\in S_N^0$,
\begin{align}\label{ch1_Z1-eq1}
(\Delta _{w} I_h u, \Delta _w v_N)  &= (I_0 u, \Delta (\Delta _w v_N )) + \langle  {I}_{g}(\nabla u)\mathbf{n_e}\cdot  {n},\Delta _w v_N\rangle \nonumber\\[2pt] 
&\qquad-\langle I_b u, \nabla (\Delta _w v_N)\cdot  {n} \rangle  \nonumber \\[2pt]
& = (\Delta I_0 u,\Delta _w v_N)+\langle ( {I}_{g}(\nabla u)\mathbf{n_e}-\nabla(I_0 u))\cdot  {n},\Delta _w v_N\rangle \nonumber \\[2pt]
&\qquad -\langle I_bu-I_0 u, \nabla (\Delta _w v_N)\cdot  {n} \rangle \nonumber \\[2pt]
& = (v_0,\Delta  (\Delta  I_0 u))-\langle v_b, \nabla (\Delta  I_0 u )\cdot  {n}\rangle +\langle (v_{g} \mathbf {n_e})\cdot \mathbf{n}, \Delta  I_0 u \rangle  \nonumber \\[2pt]
& \qquad +\langle ( {I}_{g}(\nabla u)\mathbf{n_e}-\nabla(I_0 u))\cdot  {n},\Delta _w v_N\rangle \nonumber \\[2pt]
&\qquad-\langle I_bu-I_0 u, \nabla (\Delta _w v_N)\cdot  {n} \rangle \nonumber \\[2pt]
&= (\Delta v_0 ,\Delta I_0 u)-\langle v_b-v_0, \nabla (\Delta  I_0 u )\cdot  {n}\rangle \nonumber \\[2pt]
&\qquad+\langle (v_{g} \mathbf {n_e}-\nabla v_0)\cdot \mathbf{n}, \Delta  I_0 u \rangle \nonumber \\[2pt]
&\qquad+\langle ( {I}_{g}(\nabla u)\mathbf{n_e}-\nabla(I_0 u))\cdot  {n},\Delta _w v_N\rangle \nonumber\\[2pt]
&\qquad-\langle I_bu-I_0 u, \nabla (\Delta _w v_N)\cdot  {n} \rangle \nonumber \\[2pt]
&=\big( \Delta^2 u, v_0\big)+\big( \Delta (I_0 u-u),\Delta v_0\big) \nonumber\\[2pt] 
&\qquad+\big\langle \Delta u-\Delta I_0 u, (\nabla v_0-v_g\textbf{n}_e) \cdot \textbf{n}\big\rangle \nonumber\\[2pt]
&\qquad -\big\langle \nabla(\Delta u-\Delta I_0 u)\cdot \textbf{n},v_0-v_b\big\rangle \nonumber\\[2pt]
&\qquad+\langle ( {I}_{g}(\nabla u)\mathbf{n_e}-\nabla(I_0 u))\cdot  {n},\Delta _w v_N\rangle  \nonumber \\[2pt]
&\qquad  -\langle I_bu-I_0 u, \nabla (\Delta _w v_N)\cdot  {n} \rangle.
\end{align}
Thus, we complete the proof of (\ref{ch1_2Laplace}).

 For any $v_N\in S_N^0$, integration by parts yields
\begin{eqnarray}\label{ch1_Z2-eq1}
\left( \nabla (I_0 u),\nabla v_0\right)&=& \left(\nabla\left(I_0  u-u\right), \nabla v_0\right)+\left(\nabla u, \nabla v_0\right)\nonumber \\[2pt]
&=&\left(\nabla\left(I_0   u-u\right), \nabla v_0\right)-\left(\Delta u, v_0\right)+\left\langle\nabla u \cdot \mathbf{n}, v_0\right\rangle \nonumber \\[2pt]
&=&\left(\nabla\left(I_0   u-u\right), \nabla v_0\right)-\left(\Delta u, v_0\right)\nonumber \\[2pt]
&& \, +\,\sum_{T \in \mathcal{T}_N}\left\langle\nabla u \cdot \mathbf{n}, v_0-v_b\right\rangle_{\partial T},
\end{eqnarray}
where we have used $\sum_{T \in \mathcal{T}_N}\left\langle\nabla u \cdot \mathbf{n}, v_b\right\rangle_{\partial T}=0$ due to the boundary condition $\left.v_b\right|_{\partial \Omega}=0$. Using (\ref{ch1_weakgradient}), integration by parts, and (\ref{ch1_Z2-eq1}) we have
$$
\begin{aligned}
\big( \nabla_w (I_h u),\nabla_w v_N)&=-(I_0 u, \nabla \cdot \nabla_w v_N\big)+\big\langle I_b u, \nabla_w v_N\cdot \mathbf{n}\big\rangle\\[2pt]
&=\big( \nabla (I_0 u), \nabla_w v_N \big)+\big\langle I_b u-I_0 u, \nabla_w v_N\cdot \mathbf{n}\big\rangle\\[2pt]
&=-\big( \nabla \cdot (\nabla I_0 u),v_0\big)+ \left\langle \nabla I_0 u  \cdot \mathbf{n},v_b\right\rangle +\big\langle I_b u-I_0 u, \nabla_w v_N\cdot \mathbf{n}\big\rangle\\[2pt]
& =\left(\nabla\left(I_0   u\right), \nabla v_0\right)+ \left\langle \nabla I_0 u  \cdot \mathbf{n},v_b-v_0\right\rangle +\big\langle I_b u-I_0 u, \nabla_w v_N\cdot \mathbf{n}\big\rangle \\[2pt]
& =-\left(\Delta u, v_0\right)+\left(\nabla\left(I_0   u-u\right), \nabla v_0\right) +\left\langle\left(\nabla I_0   u-\nabla u\right) \cdot \mathbf{n}, v_b-v_0\right\rangle \\
&\qquad+\big\langle I_b u-I_0 u, \nabla_w v_N\cdot \mathbf{n}\big\rangle .
\end{aligned}
$$
This completes the proof of (\ref{ch1_Laplace}). Finally, the equation (\ref{ch1_reaction})  holds trivially.
\end{proof}

From (\ref{ch1_bilinear}) and Lemma \ref{errorequation},  the error equation for the WG-FEM can be obtained in the next lemma.
\begin{lemma}
    Assume that $I_h u$ is the interpolation of the solution $u$ of \eqref{ch1_md_problem} and $u_N$ is the numerical solution computed by \eqref{ch1_bilinear}. Then, for any $v_N=\{v_0,v_b,v_g\textbf{n}_e\}\in S_N^0$, there holds:
    \begin{equation}\label{ch1_err_eqn}
        \textsl{a}(I_h u-u_N,v_N)=Z(u,v_N)+s(I_h u,v_N),
    \end{equation}
    where $Z(u,,v_N)=Z_1(u,v_N)+Z_2(u,v_N)+Z_3(u,v_N)$ and are defined as in \eqref{ch1_2Laplace}, \eqref{ch1_Laplace} and \eqref{ch1_reaction} respectively.
\end{lemma}

%%%%%%%%%%%%%%%%%%%%%%%%%%%%%%%%%%%%%%%%%%%%%%%%%%%%%%%%%%%%%%%%%%%%%%%%%%%%%%%%%%%%%%%%%%%%%%
\subsection{Error estimates in $H^{2}-$equivalent norm}
\begin{theorem}\label{H2result}
Let $I_h u$ be the interpolation of the solution $u$  of the problem (\ref{ch1_md_problem}) and  $u_{N} \in S_{N}^0$ be the Weak Galerkin finite element solution computed by the scheme \eqref{ch1_bilinear}. Let $\sigma\strut_T$ and $\rho\strut_T$ be defined by \eqref{ch1_sigma} and \eqref{ch1_rho}.
Then there exists a constant $C$ such that:
\begin{equation}\label{ch1_discrete_error}
    \||I_h u-u_N|\|\strut_{M} \leq C\Big( \varepsilon^{1/2}  N^{-(s-2)}\ln ^{s-3/2}N + N^{-(s-1)}+N^{-\lambda}\Big).
\end{equation}
\end{theorem}
\begin{proof}
 Let $\nu:= I_h u-u_N$. Recalling the definition of the norm $ \||\cdot|\|\strut_M$ and substituting $v_N=\nu$ in the equation \eqref{ch1_err_eqn}, we obtain
\begin{eqnarray}\label{ch1_errnorm}
   C \||\nu |\|\strut_{M}^{2} \leq \textsl{a}(\nu,\nu )=Z(u,\nu)+s(I_h u,\nu).
\end{eqnarray}
Next, we will estimate each term in $Z(u,\nu)$ separately and $s(I_h u,\nu)$ as follows:
Recall that %  $Z_1(u,\nu)$
\begin{eqnarray*}
    Z_1(u,\nu) &=& \varepsilon^2\big( \Delta (I_0 u - u),\Delta \nu_0\big)+\varepsilon^2\big\langle \Delta u-\Delta I_0 u, (\nabla \nu_0-v_g\textbf{n}_e) \cdot \textbf{n}\big\rangle \\ [2pt]
    &&-\varepsilon^2 \big\langle \nabla(\Delta u-\Delta I_0 u)\cdot \textbf{n},\nu_0-\nu_b\big\rangle   +\varepsilon^2\langle ( {I}_{g}(\nabla u)\mathbf{n_e}-\nabla(I_0 u))\cdot  {n},\Delta _w \nu\rangle  \\[2pt]
    &&  -\varepsilon^2\langle I_bu-I_0 u, \nabla (\Delta _w \nu)\cdot  {n} \rangle\\[2pt]
     &=:& \sum_{j=1}^5 Z_1^j(u,\nu)\,.
\end{eqnarray*}
From the Cauchy-Schwarz inequality and Lemma \ref{interpolationlemma}, one has
\begin{eqnarray}
    \vert Z_1^1(u,\nu)\vert &\leq & C \varepsilon \Vert  \Delta (u-I_0 u)\Vert  \varepsilon \Vert \Delta \nu_0 \Vert \nonumber \\ [4pt]
    &\leq & C \left(\varepsilon^{\frac{1}{2}}\left(N^{-1} \ln N\right)^{s-2}+\min \left\{\varepsilon^{\frac{1}{2}} N^2, N^{\frac{1}{2}}\right\} N^{-\lambda}\right) \Vert| \nu \Vert |\strut_{M}. \nonumber \\ \label{ch1_Z11}
\end{eqnarray}
It follows from the Cauchy-Schwarz, the H\"older inequalities, and Lemma \ref{edgeinterpolation} that
\begin{eqnarray}
\vert Z_1^2(u,\nu)\vert &\leq & \varepsilon^2 \sum_{T \in \mathcal{T}_{N}} \Vert   \Delta (u - I_{0} u)\Vert_{L^\infty(\partial T)}\Vert \nabla \nu_{0}\cdot \textbf{n}_{e}- v_{g}\Vert _{L^1(\partial T)} \nonumber \\ [2pt]
&\leq & \varepsilon^2 \sum_{T \in \mathcal{T}_{N}} \Vert   \Delta (u - I_{0} u)\Vert_{L^\infty(\partial T)} \vert \partial T\vert^{1/2}  \Vert \nabla \nu_{0}\cdot \textbf{n}_{e}- v_{g}\Vert_{L^2(\partial T)} \nonumber \\ [2pt]
&\leq & \varepsilon \Big( \sum_{T \subset \Omega_0} (  \min\{\varepsilon^{1/2}N^2, N^{1/2}\} N^{-\lambda})^2  \vert \partial T\vert \rho_T^{-1}\Big)^{1/2} \times \nonumber \\ 
&& \Big( \sum_{T \subset \Omega_0} \rho_T \Vert \nabla \nu_{0}\cdot \textbf{n}_{e}- v_{g}\Vert _{L^2(\partial T)}^2\Big)^{1/2} + \nonumber \\ [2pt]
&& \varepsilon \Big( \sum_{T \subset \Omega\setminus \Omega_0} (\varepsilon^{1/2}  (N^{-1}\ln N)^{s-2})^{2}  \vert \partial T\vert \rho_T^{-1}\Big)^{1/2} \times \nonumber \\
&&  \Big( \sum_{T \subset \Omega\setminus \Omega_0} \rho_T \Vert \nabla \nu_{0}\cdot \textbf{n}_{e}- v_{g}\Vert _{L^2(\partial T)}^2\Big)^{1/2}  \nonumber \\
&\leq & C \Big( \varepsilon^{1/2}  N^{-(s-2)}\ln ^{s-3/2}N+ \min\{\varepsilon^{1/2}N^{2}, N^{1/2}\} N^{-\lambda}\Big) \Vert| \nu \Vert |\strut_{M}. \nonumber \\ \label{ch1_Z12}
\end{eqnarray}
Similarly, we observe that
\begin{eqnarray}
\left| Z_{1}^3(u,\nu)\right| &\leq&   \varepsilon^2 \sum_{T \in \mathcal{T}_{N}} \Vert \nabla(\Delta u - \Delta I_0 u)\Vert_{L^\infty(\partial T)}\Vert  \nu_{0} - \nu_{b} \Vert _{L^1(\partial T)} \nonumber \\ [2pt]
&\leq& \varepsilon^2 \sum_{T \in \mathcal{T}_{N}} \Vert \nabla(\Delta u - \Delta I_0 u)\Vert_{L^\infty(\partial T)} \vert \partial T\vert^{1/2}   \Vert  \nu_{0} - \nu_{b}\Vert _{L^2(\partial T)} \nonumber \\ [2pt]
&\leq&  \Big( \sum_{T \subset \Omega_0} (  \min\{\varepsilon^{1/2}N^3,N\} N^{-\lambda})^2  \vert \partial T\vert \sigma_T^{-1}\Big)^{1/2}  \times \nonumber \\ [2pt]
&& \Big( \sum_{T \subset \Omega_0}\sigma_T \Vert  \nu_{0} - \nu_{b} \Vert _{L^2(\partial T)}^2\Big)^{1/2} + \nonumber \\ [2pt]
&&  \Big( \sum_{T \subset \Omega\setminus \Omega_0} (\varepsilon^{1/2}  (N^{-1}\ln N)^{s-3})^{2}  \vert \partial T\vert \sigma_T^{-1}\Big)^{1/2} \times \nonumber \\ [2pt]
&& \Big( \sum_{T \subset \Omega\setminus \Omega_0}\sigma_T \Vert\nu_{0} - \nu_{b}\Vert _{L^2(\partial T)}^2\Big)^{1/2} \nonumber \\ [2pt]
&\leq& C\Big( \varepsilon^{1/2}  N^{-(s-2)}\ln ^{s-3/2}N+ \min\{\varepsilon^{1/2}N^3,N\} N^{-\lambda}\Big) \Vert| \nu \Vert |\strut_{M}.  \nonumber \\ \label{ch1_Z13}
\end{eqnarray}
Now, we shall estimate $Z_1^4(u,\nu)$. It follows from the Cauchy-Schwarz, Lemma \ref{trace} and Lemma \ref{edgeinterpolation} that
\begin{eqnarray}\label{ch1_Z14}
\left|Z_1^4\right| &=& \varepsilon^2\left|\sum_{T \in \mathcal{T}_N} \sum_{e \subset   T}\left\langle ( {I}_{g}(\nabla u)\mathbf{n_e}-\nabla(I_0 u))\cdot  {n}, \Delta_w \nu \right\rangle_e\right| \nonumber \\ [2pt]
& \leq & \varepsilon^2 \sum_{T \in \mathcal{T}_N} \sum_{e \subset \partial T}\left\|  {I}_{g}(\nabla u)\mathbf{n_e}-\nabla(I_0 u)u\right\|_{L^{\infty}(e)}\left\|\Delta_w  \nu\right\|_{L^1(e)} \nonumber \\ [2pt]
& \leq & C \varepsilon \sum_{T \in \mathcal{T}_N}\left(\left\| {I}_{g}(\nabla u)\mathbf{n_e}-\nabla u\right\|_{L^{\infty}(\partial T)}+\left\|\nabla u- \nabla(I_0 u)\right\|_{L^{\infty}(\partial T)}\right) \times \nonumber \\ [2pt]
&& \frac{|e|}{|T|^{1 / 2}} \cdot \varepsilon\left\|\Delta_w  \nu\right\|_{L^2(T)} \nonumber \\ [2pt]
& \leq & C\Big(  N^{-(s-1)}+C  N^{-\lambda} + \varepsilon^{1/2} N^{-(s-1)} \ln ^{s-1} N\Big) \||\nu|\|\strut_{M} \, \cdot
\end{eqnarray}

Finally, we estimate the last term $Z_1^5(u,\nu)$. 

From {the H\"older inequalities, and Lemma \ref{trace}, we have for $T\in \mathcal{T}_N$
\[\left\|\nabla(\Delta_w  \nu)\right\|_{L^1(e)}\leq C\vert e\vert ^{1/2}\left\|\nabla(\Delta_w  \nu)\right\|_{L^2(e)}\leq C \kappa\dfrac{|e|}{|T|^{1/2}} \left\|\nabla(\Delta_w  \nu)\right\|_{L^2(T)}. \]}
{Again, the H\"older inequalities, an anisotropic inverse inequality  \cite[Corollary A2]{franz2014c0}  and Lemma \ref{edgeinterpolation} give}
\begin{eqnarray}\label{ch1_Z15}
\left|Z_1^5\right| &=& \varepsilon^2\left|\sum_{T \in \mathcal{T}_N} \sum_{e \subset   T}\left\langle ( {I}_{b}u-I_0 u)\cdot  {n}, \nabla(\Delta_w \nu )\right\rangle_e\right| \nonumber \\ [2pt]
& \leq & \varepsilon^2 \sum_{T \in \mathcal{T}_N} \sum_{e \subset \partial T}\left\|  {I}_{b}u-I_0 u\right\|_{{{L^\infty}}(e)}\left\|\nabla(\Delta_w  \nu)\right\|_{{L^1}(e)} \nonumber\\ [2pt]
& \leq &  C \varepsilon \sum_{T \in \mathcal{T}_N}\left(\left\| {I}_{b}u-u\right\|_{{L^\infty}(\partial T)}+\left\| u- I_0 u\right\|_{{L^\infty}(\partial T)}\right) \times \nonumber \\ [2pt]
&& { \kappa \dfrac{|e|}{|T|^{1/2}}h_{T}^{-1}} \cdot \varepsilon\left\|\Delta_w  \nu\right\|_{L^2(T)} \nonumber \\
& \leq & \varepsilon \sum_{T \subset \Omega_0} (  N^{-s}+\varepsilon N^{-\lambda})   N {  \varepsilon\left\|\Delta_w  \nu\right\|_{L^2(T)}} + \nonumber  \\ [2pt]
&& \varepsilon \sum_{T \subset \Omega\setminus \Omega_0} \varepsilon^{3/2}  (N^{-1}\ln N)^{s}  (\varepsilon N^{-1} \ln N)^{-1}   {  \varepsilon\left\|\Delta_w  \nu\right\|_{L^2(T)}} \nonumber \\ [2pt]
&\leq & C\Big( {N^{-s}+\varepsilon N^{-\lambda}} +\varepsilon^{1/2}  (N^{-1}\ln N)^{s-1} \Big) \Vert| \nu \Vert |\strut_{M} .
\end{eqnarray}
Now we will bound $Z_{2}(u,\nu)$ and $Z_{3}(u,\nu)$ as follows.
First
\begin{eqnarray*}
    Z_{2}(u,\nu) &=& \left(\nabla\left(u-I_0 u\right), \nabla \nu_0\right)+ \left\langle\left(\nabla u-\nabla I_0 u\right) \cdot \mathbf{n}, \nu_b-\nu_0\right\rangle \\ [2pt]
    && +\big\langle I_0 u-I_b u, \nabla_w \nu_N\cdot \mathbf{n}\big\rangle \\ [2pt]
    &=:& \sum_{j=1}^3 Z_2^j(u,\nu).
\end{eqnarray*}
From the Cauchy-Schwarz inequality, and Lemma \ref{interpolationlemma}, we have
\begin{eqnarray}
    |Z_2^1(u,\nu)| &\leq & \|\nabla\left(u-I_0 u\right)\|\| \nabla \nu_0\| \nonumber \\ [2pt]
     &\leq & C\left(\varepsilon^{1/2}\left(N^{-1} \ln N\right)^{s-1}+N^{-(s-1)}+N^{-\lambda}\right)^2\||\nu|\|\strut_{M}. \nonumber \\ \label{Z21}
\end{eqnarray}
Using the Cauchy-Schwarz inequality, H\"older inequalities and Lemma \ref{edgeinterpolation} we can bound
\begin{eqnarray}
    |Z_2^2(u,\nu)| &\leq & \sum_{T \in \mathcal{T}_{N}}^{}\|\nabla\left(u-I_0 u\right)\|_{L^{\infty}(\partial T)}\|\nu_b - \nu_0\|_{L^{1}(\partial T)} \nonumber \\ [2pt]
    &\leq & \sum_{T \in \mathcal{T}_{N}}^{}\|\nabla\left(u-I_0 u\right)\|_{L^{\infty}(\partial T)}\,|\partial T|^{1/2}\|\nu_b - \nu_0\|_{L^{2}(\partial T)}  \nonumber \\ [2pt]
    &\leq& \left(\sum_{T \subset \Omega_0}C (N^{-(s-1)}+  N^{-\lambda})^{2}|\partial T|\sigma_T^{-1}\right)^{1/2} \times \nonumber \\ [2pt]
    && \left(\sum_{T \subset \Omega_0} \sigma_T\|\nu_b-\nu_0\|^{2}_{L^{2}(\partial T)}\right)^{1/2} + \nonumber \\ [2pt]
    && \left(\sum_{T \subset \Omega \setminus \Omega_0}^{}C (\varepsilon^{1/2} N^{-(s-1)} \ln ^{s-1} N)^{2}|\partial T| \sigma_T^{-1}\right)^{1/2} \times \nonumber \\ [2pt]
    && \left(\sum_{T \subset \Omega \setminus \Omega_0} \sigma_T \|\nu_b-\nu_0\|^{2}_{L^{2}(\partial T)}\right)^{1/2} \nonumber \\ [5pt]
    & \leq & C\left(\varepsilon^{1/2} N^{-s} \ln ^{s+1/2} N+N^{-(s-1)}+  N^{-\lambda}\right)\||\nu|\|\strut_{M}.
\end{eqnarray}
Now using triangle inequality, H\"older inequalities and Lemma \ref{trace}
\begin{eqnarray}
|Z_2^3(u,\nu)|  &\leq& \left|\sum_{T \in \mathcal{T}_{N}}^{} \sum_{e \subset \partial T}^{}\big\langle I_0 u-I_b u, \nabla_w \nu_N\cdot \mathbf{n}\big\rangle_{e}\right| \nonumber \\ [2pt]
 & \leq & \sum_{T \in \mathcal{T}_{N}}^{} \sum_{e \subset \partial T}^{}\|I_0 u-I_b u\|_{L^{\infty}(e)}\|\nabla_w \nu_N\|_{L^{1}(e)} \nonumber \\ [2pt]
& \leq & C\sum_{T \in \mathcal{T}_{N}}^{}\left(\|I_0 u - u\|_{L^{\infty}(\partial T)} + \|u - I_b u\|_{L^{\infty}(\partial T)}\right)\frac{|e|}{|T|^{1 / 2}} \|\nabla_w \nu_N\|_{L^{2}(T)} \nonumber \\ [2pt] 
 & \leq & C\left((N^{-s}+\varepsilon N^{-\lambda})N^{1/2} + \varepsilon^{1/2} N^{-s} \ln^{s-1/2}N\right)\||\nu|\|\strut_{M} \, \cdot
\end{eqnarray}
with this we have done the evaluation for finding the bound for $Z_2(u, \nu)$ and now we need to bound $Z_3(u,\nu)$.
We know that
\begin{eqnarray}\label{ch1_Z3}
        %&Z_3(u,\nu) = (a(u-I_0 u),\nu_0)\\ \nonumber
        |Z_3(u,\nu)| \, &\leq & C\|u-I_0 u\|_{L^{2}(\Omega)} \|a\nu_0\|_{L^{2}(\Omega)} \nonumber \\ 
        &\leq & C\left(\varepsilon^{\frac{3}{2}}\left(N^{-1} \ln N\right)^s+N^{-s}+\varepsilon N^{-\lambda}\right)^2\|a\nu_0\|\strut_{L^{2}(\Omega)} \nonumber \\
        &\leq & C\left(\varepsilon^{\frac{3}{2}}\left(N^{-1} \ln N\right)^s+N^{-s}+\varepsilon N^{-\lambda}\right)^2\||\nu|\|\strut_{M}.
\end{eqnarray}
using the Cauchy-Schwarz inequality and Lemma \ref{interpolationlemma}.
Now recall that our stabilizer term is given by
\begin{equation*}
    s(I_h u,\nu) = \sum_{T \in \mathcal{T}_{N}}^{}\left[\rho_T\langle\nabla I_{0} u - I_{g}{(\nabla u \cdot \textbf{n}_e)} \,\textbf{n}_{e}, \nabla \nu_{0}-\nu_{g}\textbf{n}_{e}\rangle\strut_{\partial T} + \sigma_T\langle I_{0} u - I_{b} u, \nu_{0}-\nu_{b} \rangle\strut_{\partial T}\right].
\end{equation*}
Then,
\begin{eqnarray}\label{ch1_stabeq}
    |s(I_h u, \nu)| &\leq& \sum_{T \in \mathcal{T}_{N}}^{}\Big(|\rho_T\langle\nabla I_{0} u\cdot \textbf{n}_{e}- I_{g}{(\nabla u \cdot \textbf{n}_e) } , \nabla \nu_{0}\cdot\textbf{n}_{e}-\nu_{g}\rangle\strut_{\partial T}| \nonumber \\ [2pt]
    && \quad+ |\sigma_T\langle I_{0} u - I_{b} u, \nu_{0}-\nu_{b} \rangle\strut_{\partial T}|\Big) \nonumber \\ [2pt]
    & \leq & \sum_{T \in \mathcal{T}_{N}}^{} \rho_T\|\nabla I_{0} u\cdot \textbf{n}_{e}- I_{g} {(\nabla u \cdot \textbf{n}_e)}\|\strut_{L^{\infty}(\partial T)}\|\nabla \nu_{0}\cdot \textbf{n}_{e}-\nu_{g}\|\strut_{L^{1}(\partial T)} \nonumber \\  [2pt]
    && + \sum_{T \in \mathcal{T}_{N}}^{} \sigma_T\|I_{0} u - I_{b} u\|\strut_{L^{\infty}(\partial T)}\|\nu_{0}-\nu_{b}\|\strut_{L^{1}(\partial T)} \nonumber \\ [2pt]
    & \leq & \sum_{T \in \mathcal{T}_{N}}^{} \sigma_T^{1/2}\left(\|I_0 u - u\|\strut_{L^{\infty}(\partial T)} + \|u - I_b u\|\strut_{L^{\infty}(\partial T)}\right)|\partial T|^{1/2} \times \nonumber \\ [2pt] 
    && \sigma_T^{1/2}\|\nu_{0}-\nu_{b}\|\strut_{L^{2}(\partial T)} + \nonumber \\ [2pt]
    && \sum_{T \in \mathcal{T}_{N}}^{} \rho_T^{1/2}\left(\|\nabla I_{0} u - \nabla u\|\strut_{L^{\infty}(\partial T)} + \|\nabla u - I_g {(\nabla u \cdot \textbf{n}_e)} \textbf{n}_{e}\|\strut_{L^{\infty}(\partial T)}\right) \times \nonumber \\ [2pt]
    &&  |\partial T|^{1/2} \cdot \rho_T^{1/2}\|\nabla \nu_{0}-\nu_{g}\textbf{n}_{e}\|\strut_{L^{2}(\partial T)} \nonumber \\ [2pt]
    & \leq & C\left(N^{-s}+\varepsilon N^{-\lambda} + \varepsilon^{3/2} N^{-(s-3/2)} \ln^{s-3/2}N\right)\||\nu|\|\strut_{M} \nonumber \\ [2pt]
    &&  + C\Big(  N^{-(s-1)}+  N^{-\lambda} + \varepsilon^{1/2} N^{-(s-3/2)} \ln ^{s-3/2} N\Big) \||\nu|\|\strut_{M}\, \cdot
\end{eqnarray}
Now using the above estimates \eqref{ch1_Z11}-\eqref{ch1_stabeq} and \eqref{ch1_errnorm} we have the desired error bound. Hence,  we complete the proof.
\end{proof}

\begin{remark}
 Although the anisotropic inverse inequality given in \cite[Corollary A2]{franz2014c0} is for rectangular elements, the same inequality holds true for triangular elements if the coordinate-system condition is satisfied. For more details one can refer to \cite[Lemma 3.1]{apel1994local}.
\end{remark}

Finally, we are ready to state and prove the main result of the paper in the following theorem:
\begin{theorem}\label{maintheorem}
    Let $u$ be the solution  of the problem (\ref{ch1_md_problem}) and  $u_{N} \in S_{N}^0$ be the Weak Galerkin finite element solution computed by the scheme \eqref{ch1_bilinear}. Let $\rho\strut_T$ and $\sigma\strut_T$ be defined by \eqref{ch1_sigma} and \eqref{ch1_rho}. Assume that $\lambda \geq s$, then there exists a constant $C$ such that
\begin{eqnarray}\label{ch1_main}
    \| | u-u_N|\|\strut_{M} \,\, \leq C\Big( \varepsilon^{1/2}  N^{-(s-2)}\ln ^{s-3/2}N + N^{-(s-1)}\Big).
\end{eqnarray}
\end{theorem}
%{\bf Proof.}
\begin{proof}
We will use the triangle inequality to obtain
\begin{equation*}
|||u-u_N|||\strut_{M} \leq \,\, |||u-I_{h}u|||\strut_{M} + |||I_{h}u - u_N|||\strut_{M} \,\, ,
\end{equation*}
where by Theorem \ref{H2result} we know that $|||I_{h}u - u_N|||\strut_{M}$ satisfies \eqref{ch1_main}.
Now we need to bound $|||u-I_{h}u|||\strut_{M}$ and so by using equation \eqref{ch1_triplenorm} and anisotropic interpolation estimates in Lemma \ref{interpolationlemma} one can get the following result:
\begin{align*}
\sum_{T\in \mathcal{T}_N}\Big(\varepsilon^2\Vert \Delta (u-I_0u)\Vert\strut_T^2 &+ \Vert\nabla (u-I_0 u)\Vert\strut_T^2 + |a|\Vert (u-I_{0}u)\Vert\strut_T^2\Big)\\
&\leq C\Big( \varepsilon^{1/2}  N^{-(s-2)}\ln ^{s-3/2}N + N^{-(s-1)}\Big).
\end{align*}
Further, using Lemma \ref{edgeinterpolation}, we get
\begin{equation*}  
\sum_{\partial T \in \mathcal{T}_{N}}^{} s(u - I_h u, u - I_h u) \leq C\Big( \varepsilon^{1/2}  N^{-(s-2)}\ln ^{s-3/2}N + N^{-(s-1)}\Big), 
\end{equation*}
which completes the proof.
\end{proof}

\section{Numerical Results}\label{Numerical_Examples}

\begin{figure}[!h]
\centerline{%
	\begin{tabular}{cc}
		\resizebox*{7cm}{!}{\includegraphics{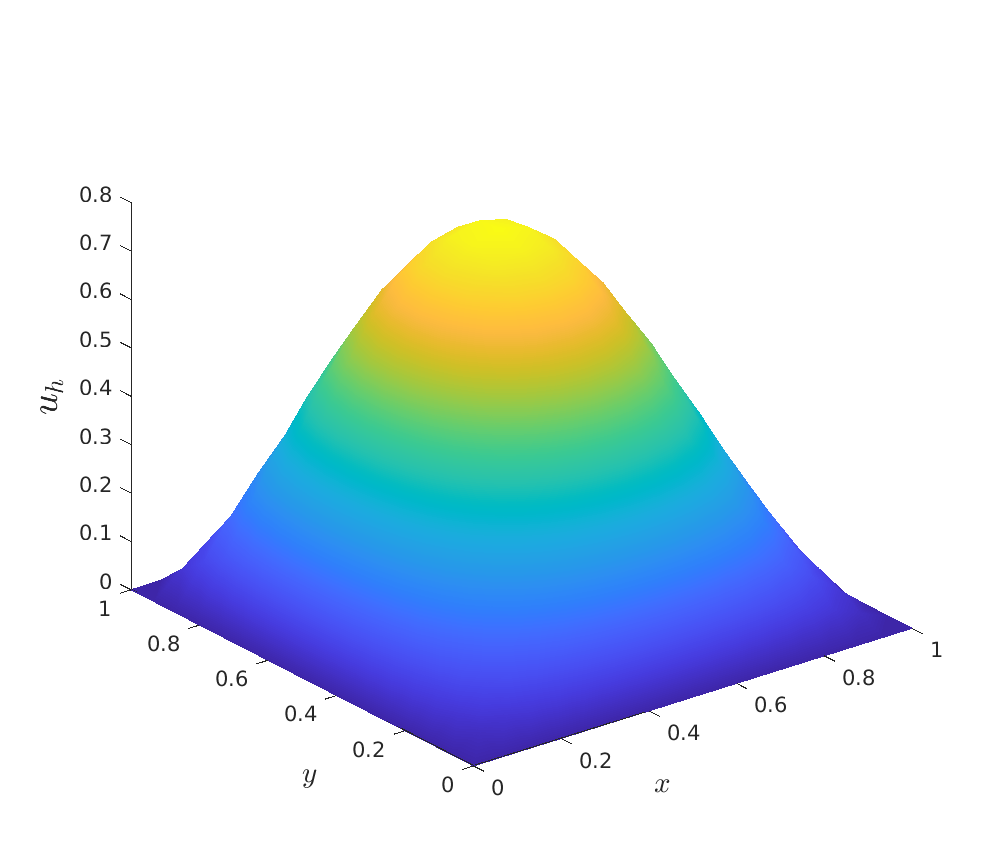}}
		&
		\resizebox*{7cm}{!}{\includegraphics{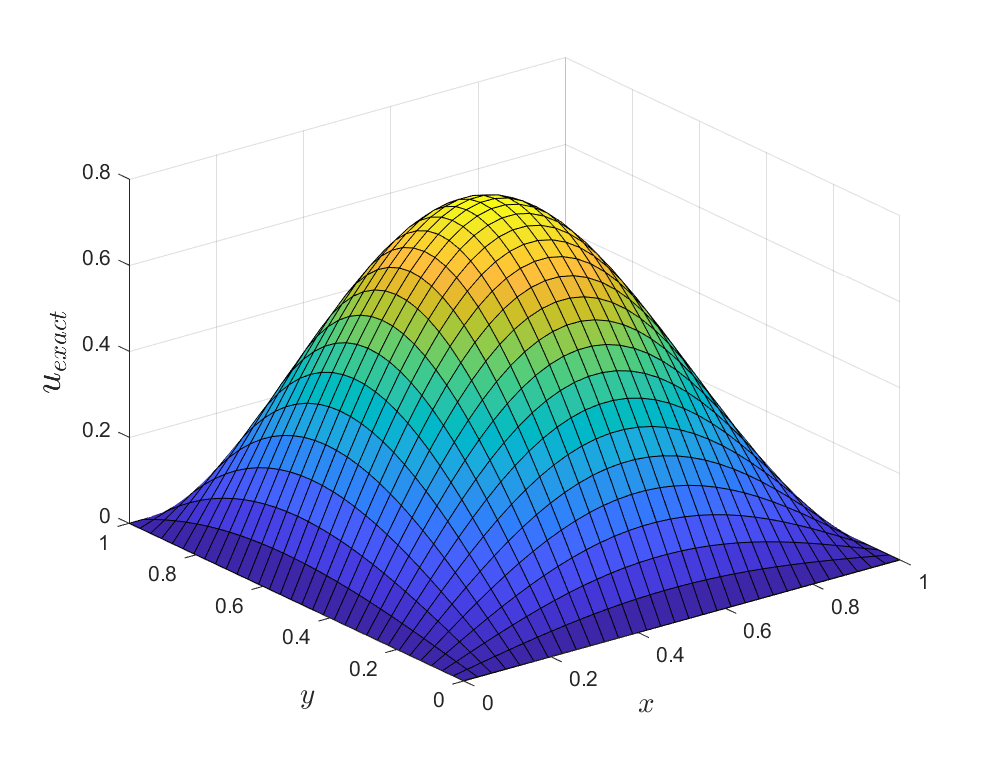}}\\
		{\it {\em (a) WG Solution.}} & {\it {\em (b) Exact Solution.}}
	\end{tabular}
} \caption{ Comparison of solution for Example \ref{ex1} for $\varepsilon^2 = 10^{-6}$ and $N=32$. }\label{ch1_fig1}
\end{figure}

In this section, we will test the proposed WG-FEM for some chosen examples and will show the obtained results. We will be using the following discrete finite element space in this section:
\begin{align*}
\tilde{S}_{N}(\Omega) = \{u = \{u_{0},u_{b},u_{g}\textbf{n}_{e}\} :\, &u_{0} \in P_{{k}}(T), u_{b} \in P_{{k}}(e),\\ \nonumber
&u_{g} \in P_{{k-1}}(e), \,\, \mbox{where}\,\, T \in \mathcal{T}_{N}, e \in \partial T \, \, \mbox{and} \, \, {k = 2,3} \}.
\end{align*}
\begin{example}\label{ex1}
Consider the problem \eqref{ch1_md_problem} with $a = 0$ in a unit square $\Omega = (0,1)^{2}$. We choose $\mathrm{g}$ in a way that the exact solution is
\begin{eqnarray*}
    u(x,y) &=& \left[\sin(\pi x) + \frac{\pi \varepsilon}{(1-e^{-1/\varepsilon})}(e^{-x/\varepsilon}+e^{(x-1)/\varepsilon}-1-e^{-1/\varepsilon})\right] \times \\ [2pt]
    && \Big[2 y(1-y^{2}) + \varepsilon\Big(\delta_{1}\delta_{2}(1-2y)-\frac{3\delta_{2}}{\delta_{1}}+(\frac{3}{\delta_{1}}-\delta_{2})e^{-y/\varepsilon}\\ [2pt]
    && \qquad \qquad \qquad \qquad +(\frac{3}{\delta_{1}}+\delta_{2})e^{(y-1)/\varepsilon}\Big)\Big]
\end{eqnarray*}
where $\delta_{1} = 1 - e^{-1/\varepsilon}, \delta_{2} = 1 + e^{-1/\varepsilon}  \,\, \mbox{and}\,\, \delta_{3} = \dfrac{1}{\delta_{2}-2\varepsilon\delta_{1}}$. The error in  $\||\cdot|\|$ norm and order of convergence is given in Table \ref{ch1_table1}. Comparsion of WG solution and Exact solution is given in Figure \ref{ch1_fig1}.
\begin{table}[htbp]
\centering
\caption{Error in the energy norm and Convergence rates for Example \ref{ex1}.}\label{ch1_table1}
\begin{tabular}{llllll}
\toprule
&{$N$}&\makebox[3em]{$\varepsilon^2 = 10^{-6}$}&\makebox[3em]{order}&{$\varepsilon^2 = 10^{-10}$}&{order}\\
\midrule
&&\makebox[3em]{$\|| \nu|\|$}&\makebox[3em]{}&{$\||\nu|\|$}\\
\midrule
&$4$ &2.32507e+00& - &{4.1265e+00}& - \\
&$8$ &6.57597e-01& 1.822 &{1.23918e+00}& {1.735}\\
$k=2$&$16$ &1.73051e-01 &1.925&{3.45618e-01}&{1.842}\\
&$32$ &4.4421e-02&1.969& {8.91642e-02}&{1.954}\\
&{$64$} &{1.21005e-02}&{1.876}&{2.23968e-02}&{1.993}\\
&{$128$} & {3.78377e-03
} & {1.677} &{5.60062e-03}& {1.999} \\
\midrule
&{$4$} &{3.92034e+00}&-&{1.1445e+00}&-\\
&{$8$} &{5.44567e-01}&{2.847}&{1.56523e-01}&{2.870}\\
{$k=3$}&{$16$} &{7.53413e-02}&{2.853}&{2.10155e-02}&{2.896}\\
&{$32$} &{9.98824e-03}&{2.915}&{2.73532e-03}&{2.941}\\
&{$64$} &{1.78349e-04}&{2.846}&{3.42472e-04}&{2.997}\\
&{$128$} &{2.67346e-05}&{2.737}&{4.27983e-05}&{3.000}\\
\bottomrule
\end{tabular}
\end{table}
\end{example}

\begin{figure}[htbp]
 \centerline{%
 	\begin{tabular}{cc}
		\resizebox*{7cm}{!}{\includegraphics{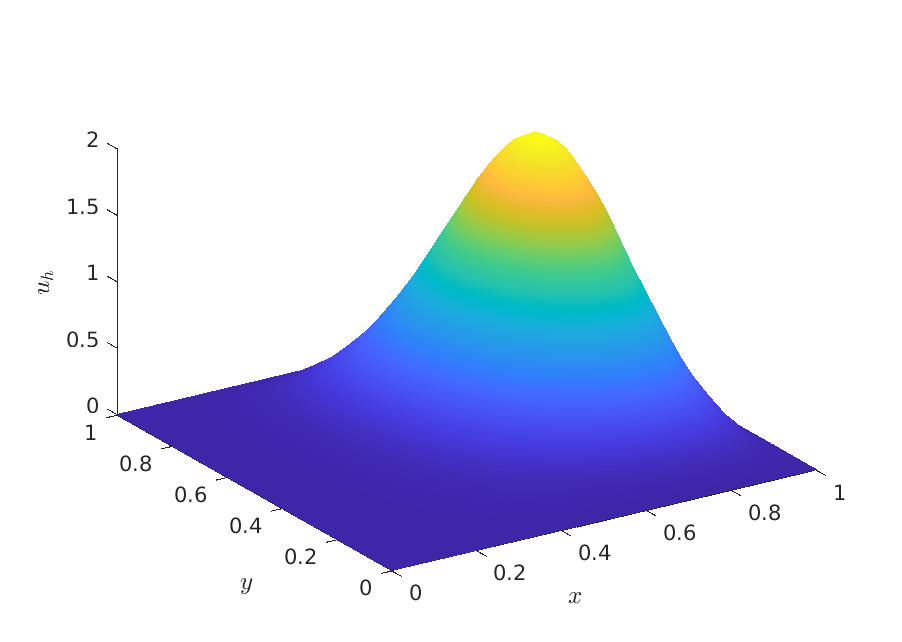}}
		&
		\resizebox*{7cm}{!}{\includegraphics{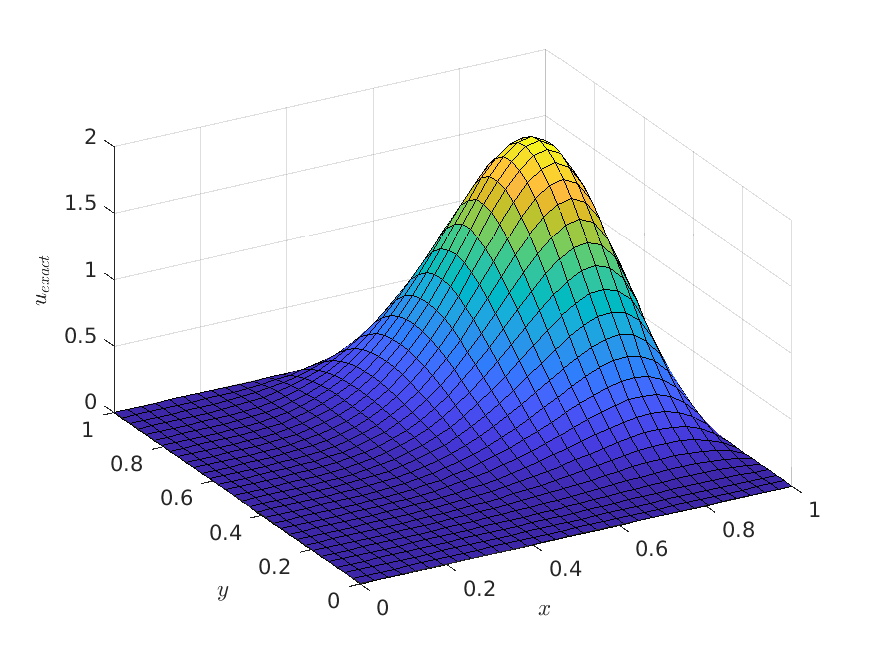}}\\
		{\it {\em (a) WG Solution.}} & {\it {\em (b) Exact Solution.}}
	\end{tabular}
} \caption{ Comparison of solution for Example \ref{ex2} for $\varepsilon^2 = 10^{-6}$ and $N=32$. }\label{ch1_fig2}
\end{figure}

{
\begin{example}\label{ex2}
In this example, we consider the problem \eqref{ch1_md_problem} with $a=x$ and choose source function in such a way that
\begin{eqnarray*}
    u(x,y) &=& 250\Big[\varepsilon(e^{-x/\varepsilon} + e^{-y/\varepsilon}) - x^2y\Big]\times\Big[\varepsilon(e^{-(1-x)/\varepsilon} + e^{-(1-y)/\varepsilon}) - x^2y\Big] \times\\
    &&  \,\, xy(1-x)(1-y)
\end{eqnarray*}
is the exact solution. The error in  $\||\cdot|\|$ norm and order of convergence is given in Table \ref{ch1_table2}. Surface plots of WG solution and Exact solution are presented in Figure \ref{ch1_fig2}.
\end{example}
\begin{remark}
    From the two examples considered, we can see that for the smaller perturbation parameter values like $\varepsilon = 10^{-5}$, the energy norm error is dominated by $N^{-(s-1)}$ which supports our theoretical result stated in Theorem \ref{maintheorem}.
\end{remark}
}

\begin{table}[htbp]
\centering
\caption{Error in the energy norm and Convergence rates for Example \ref{ex2}.}\label{ch1_table2}
\begin{tabular}{llllll}
\toprule
&{{$N$}}&\makebox[3em]{{$\varepsilon^2 = 10^{-6}$}}&\makebox[3em]{{order}}&{{$\varepsilon^2 = 10^{-10}$}}&{{order}}\\
\midrule
&&\makebox[3em]{{$\|| \nu|\|$}}&\makebox[3em]{}&{{$\||\nu|\|$}}\\
\midrule
&{$4$} &{6.68665e+00}& -&{2.19564e+00}& - \\
&{$8$} &{1.32888e+00}& {2.331} &{5.50095e-01}& {1.996}\\
{$k=2$}&{$16$} &{3.32249e-01} &{1.999}&{1.35582e-01}&{2.020}\\
&{$32$} &{8.35620e-02}&{1.991}& {3.43142e-02}&{1.982}\\
&{$64$} &{2.09243e-02}&{1.997}&{8.64901e-03}&{1.988}\\
&{$128$} & {6.15928e-03} & {1.764} &{2.16952e-03}& {1.995} \\
\midrule
&{$4$} &{1.43416e+00}&-&{8.51156e-01}&-\\
&{$8$} &{1.71222e-01}&{3.006}&{1.09358e-01}&{2.960}\\
{$k=3$}&{$16$} &{2.28266e-02}&{2.907}&{1.36778e-02}&{2.999}\\
&{$32$} &{3.21095e-03}&{2.829}&{ 1.72992e-03}&{2.983}\\
&{$64$} &{4.79175e-04}&{2.744}&{2.19005e-04}&{2.981}\\
&{$128$} &{7.71742e-05
}&{2.634}&{ 2.75634e-05}&{2.990}\\
\bottomrule
\end{tabular}
\end{table}

%%%%%%%%%%%%%%%%%%%%%%%%%%%%%%%%%%%%%%%%%%%%%%%%%%%%%%%%%%%%%%%%%%%%%%%%%%
\section{Conclusions}\label{Conclusions}
In this article, we presented a WG finite element method for solving the singularly perturbed biharmonic elliptic problem with homogeneous clamped boundary conditions. By choosing a suitable weak function, weak gradient and weak Laplacian space, we constructed our scheme and anisotropic error estimates in $H^{2}-$equivalent norm were also derived. Numerical examples are presented showing uniform convergence of the method.

\goodbreak\noindent
\section*{Data Availability}
Data sharing is not applicable to this paper because no datasets were created or
examined during the current study.

\section*{Declaration of competing interest}
The authors declare that they have no known financial interests or personal relationships that could have appeared to influence the work reported in this paper.

\section*{Acknowledgments}
The first author would like to express the thanks to Indian Institute of Technology Guwahati, India for funding of this project. \\

The authors wish to acknowledge the anonymous referees for carefully reading the manuscript and providing their valuable comments and suggestions, which really helped to improve the presentation.

%\section*{References}	

\end{document}